\documentclass[letterpaper,10pt,conference]{ieeeconf}  
\usepackage{graphicx}

\usepackage{amsthm}
\usepackage[cmex10]{amsmath}
\usepackage{amsfonts,amssymb}
\usepackage{times}
\usepackage{bm}
\usepackage{color}
\usepackage{graphicx}
\usepackage{subcaption}
\usepackage{amsfonts}
\usepackage{cite}
\usepackage{mathtools}
\IEEEoverridecommandlockouts
\providecommand{\U}[1]{\protect\rule{.1in}{.1in}}

\theoremstyle{definition}

\title{\LARGE \bf
Optimal Event-Driven Multi-Agent Persistent Monitoring of a Finite Set of Targets
}
\author{ Nan Zhou$^1$, Xi Yu$^2$, Sean B. Andersson$^{1,2}$, and Christos G. Cassandras$^{1,3}$\\
\small $^{1}$Division of Systems Engineering, $^{2}$Department of Mechanical Engineering , $^3$Department of Electrical and Computer Engineering\\
Boston University, Boston, MA 02215, USA\\
E-mail:\tt \{nanzhou,xyu,sanderss,cgc\}@bu.edu 
\thanks{* The work of Cassandras and Zhou is supported in part by NSF under grants CNS-1239021, ECCS-1509084, and IIP-1430145, by AFOSR under grant FA9550-15-1-0471, and by ONR under grant N00014-09-1-1051. The work of Andersson and Yu is supported in part by the NSF through grant ECCS-1509084.}
}

\begin{document}

\maketitle

\begin{abstract}
We consider the problem of controlling the movement of multiple cooperating
agents so as to minimize an uncertainty metric associated with a finite number
of targets. In a one-dimensional mission space, we adopt an optimal control
framework and show that the solution is reduced to a simpler parametric
optimization problem: determining a sequence of locations where each agent may
dwell for a finite amount of time and then switch direction. This amounts to a
hybrid system which we analyze using Infinitesimal Perturbation Analysis (IPA)
to obtain a complete on-line solution through an event-driven gradient-based
algorithm which is also robust with respect to the uncertainty model used. The
resulting controller depends on observing the events required to excite the
gradient-based algorithm, which cannot be guaranteed. We solve this problem by
proposing a new metric for the objective function which creates a potential
field guaranteeing that gradient values are non-zero. This approach is
compared to an alternative graph-based task scheduling algorithm for
determining an optimal sequence of target visits. Simulation examples are
included to demonstrate the proposed methods.

\end{abstract}




\section{Introduction}

\label{sec:intro} Systems consisting of cooperating mobile agents are often
used to perform tasks such as coverage control
\cite{zhong2011distributed,sun2015optimal}, surveillance, and environmental
sampling. The \emph{persistent monitoring}\ problem arises when agents must
monitor a dynamically changing environment which cannot be fully covered by a
stationary team of agents. Thus, persistent monitoring differs from
traditional coverage tasks due to the perpetual need to cover a changing
environment \cite{cassandras2013optimal,lin2013optimal}. A result of this
exploration process is the eventual discovery of various \textquotedblleft
points of interest\textquotedblright\ which, once detected, become
\textquotedblleft targets\textquotedblright\ or \textquotedblleft data
sources\textquotedblright\ which need to be monitored. This setting arises in
multiple application domains ranging from surveillance, environmental
monitoring, and energy management \cite{michael2011persistent,
smith2011persistent} down to nano-scale systems tasked to track fluorescent or
magnetic particles for the study of dynamic processes in bio-molecular systems
and in nano-medical research \cite{Shen:2011kj,cromer2011tracking}. In
contrast to \cite{cassandras2013optimal,lin2013optimal} where \emph{every}
point in a mission space must be monitored, the problem we address here
involves a \emph{finite number} of targets (typically larger than the number
of agents) which the agents must cooperatively monitor through periodic visits.

Each target may be viewed as a dynamic system in itself whose state is
observed by agents equipped with sensing capabilities (e.g., cameras) and
which are normally dependent upon their physical distance from the target. The
objective of cooperative persistent monitoring in this case is to minimize an
overall measure of uncertainty about the target states. This may be
accomplished by assigning agents to specific targets or by designing motion
trajectories through which agents reduce the uncertainty related to a target
by periodically visiting it (and possibly remaining at the target for a finite
amount of time). Viewed as an optimization problem, the goal is to jointly
minimize some cost function that captures the desired features of the
monitoring problem \cite{horling2004survey}. As long as the numbers of agents
and targets is small, it is possible to identify sequences that yield a
globally optimal solution; in general, however, this is a computationally
complex procedure which does not scale well \cite{yu2015persistent}.

Rather than viewing this problem as a scheduling task which eventually falls
within the class of traveling salesman or vehicle routing problems
\cite{stump2011multi}, in this paper we follow earlier work in
\cite{lin2013optimal} and introduce an optimal control framework whose
objective is to control the movement of agents so as to collect information
from targets (within agent sensing ranges) and ultimately minimize an average
metric of uncertainty over all targets. An important difference between the
persistent monitoring problem in previous work \cite{cassandras2013optimal}
and the current setting is that there is now a finite number of targets that
agents need to monitor as opposed to every point in the mission space. In a
one-dimensional mission space, we show that the optimal control problem can be
reduced to a parametric optimization problem. In particular, every optimal
agent trajectory is characterized by a finite number of points where the agent
switches direction and by a dwelling time at each such point. As a result, the
behavior of agents under optimal control is described by a hybrid system. This
allows us to make use of Infinitesimal Perturbation Analysis (IPA)
\cite{cassandras2010perturbation,wardi2010unified} to determine on-line the
gradient of the objective function with respect to these parameters and to
obtain a (possibly local) optimal trajectory. Our approach exploits an
inherent property of IPA which allows virtually arbitrary stochastic effects
in modeling target uncertainty. Moreover, IPA's event-driven nature renders it
\emph{scalable} in the number of events in the system and not its state space. 

A potential drawback of event-driven control methods is that they obviously
depend on the events which \textquotedblleft excite\textquotedblright\ the
controller being observable. However, this is not guaranteed under every
feasible control: it is possible that no such events are excited, in which
case the controller may be useless. The crucial events in our case are
\textquotedblleft target visits\textquotedblright\ and it is possible that
such events may never occur for a large number of feasible agent trajectories
which IPA uses to estimate a gradient on-line. At the heart of this problem is
the fact that the objective function we define for a persistent monitoring
problem has a non-zero cost metric associated with only a subset of the
mission space centered around targets, while all other points have zero cost,
since they are not \textquotedblleft points of interest\textquotedblright.
This lack of event excitation is a serious problem in many trajectory planning
and optimization tasks \cite{Schwager2009decentralized, cao2011maintaining, oh2014formation}. In this paper we solve this
problem using a new cost metric introduced in \cite{YasamanKhazaeni} which
creates a potential field guaranteeing that gradient values are generally
non-zero throughout the mission space and ensures that all events are
ultimately excited.

The rest of the paper is organized as follows. Section II formulates the
optimal control problem and Section III presents a Hamiltonian analysis which
characterizes the optimal solution in terms of two parameter vectors
specifying switching points and associated dwelling times. In Section IV we
provide a complete solution obtained through event-driven IPA gradient
estimation, and solve the problem of potential lack of event excitation
through a modified cost metric. Section V presents our graph-based scheduling
approach and Section VI includes several simulation results.

\section{Persistent Monitoring Problem Formulation}

We consider $N$ mobile agents moving in a one dimensional mission space
$[0,L]\subset\mathbb{R}$. Let the position of the agents at time $t$ be
$s_{j}(t)\in\left[  0,L\right]  $, $j=1,\ldots,N$, following the dynamics:%
\begin{equation}
\dot{s}_{j}(t)=u_{j}(t)\label{eq:multiDynOfS}%
\end{equation}
i.e., we assume that the agent can control its direction and speed. Without
loss of generality, after proper rescaling, we further assume that the speed
is constrained by $\left\vert u_{j}\left(  t\right)  \right\vert \leq1$,
$j=1,\ldots,N$. As will become clear, the agent dynamics in
(\ref{eq:multiDynOfS}) can be replaced by a more general model of the form
$\dot{s}_{j}(t)=g_{j}(s_{n})+b_{j}u_{j}(t)$ without affecting the main results
of our analysis. Finally, an additional constraint may be imposed if we assume
that the agents are initially located so that $s_{j}\left(  0\right)
<s_{j+1}\left(  0\right)  $, $j=1,\ldots,N-1$, and we wish to prevent them
from subsequently crossing each other over all $t$:%
\begin{equation}
s_{j}\left(  t\right)  -s_{j+1}\left(  t\right)  \leq 0\label{eq:multiNoCross}%
\end{equation}
The ability of an agent to sense its environment is modeled by a function
$p_{j}(x,s_{j})$ that measures the probability that an event at location
$x\in\left[  0,L\right]  $ is detected by agent $j$. We assume that
$p_{j}(x,s_{j})=1$ if $x=s_{j}$, and that $p_{j}(x,s_{j})$ is monotonically
nonincreasing in the distance $|x-s_{j}|$, thus capturing the reduced
effectiveness of a sensor over its range which we consider to be finite and
denoted by $r_{j}$. Therefore, we set $p_{j}(x,s_{j})=0$ when $|x-s_{j}%
|>r_{j}$. Although our analysis is not affected by the precise sensing model
$p_{j}(x,s_{j})$, we will limit ourselves to a linear decay model as follows:%
\begin{equation}
p_{j}(x,s_{j})=\max\{1-\dfrac{|s_{j}-x|}{r_{j}},0\}\label{eq:SensingCap_1A}%
\end{equation}
Unlike the persistent monitoring problem setting in
\cite{cassandras2013optimal}, here we consider a known finite set of targets
located at $x_{i}\in(0,L),$ $i=1,\ldots,M$ (we assume $M>N$ to avoid
uninteresting cases where there are at least as many agents as targets, in
which case every target can be assigned to at least one agent). We can then
set $p_{j}(x_i,s_{j}\left(  t\right)  )\equiv p_{ij}(s_{j}\left(  t\right)  )$
to represent the effectiveness with which agent $j$ can sense target $i$ when
located at $s_{j}\left(  t\right)  $. Accordingly, the joint probability that
$x_{i}\in\left(  0,L\right)  $ is sensed by all $N$ agents simultaneously
(assuming detection independence) is
\begin{equation}
P_{i}(\mathbf{s}(t))=1-\prod_{j=1}^{N}[1-p_{ij}(s_{j}%
(t))]\label{eq:SensingCap_NA}%
\end{equation}
where we set $\mathbf{s}(t)=[s_{1}\left(  t\right)  ,\ldots,s_{N}\left(
t\right)  ]^{\text{T}}$. Next, we define uncertainty functions $R_{i}(t)$
associated with targets $i=1,\ldots,M$, so that they have the following
properties: $(i)$ $R_{i}(t)$ increases with a prespecified rate $A_{i}$ if
$P_{i}\left(  \mathbf{s}(t)\right)  =0$ (we will later allow this to be a
random process $\{A_{i}(t)\}$), $(ii)$ $R_{i}(t)$ decreases with a fixed rate
$B_{i}$ if $P_{i}\left(  \mathbf{s}(t)\right)  =1$ and $(iii)$ $R_{i}(t)\geq0$
for all $t$. It is then natural to model uncertainty dynamics associated with
each target as follows:
\begin{equation}
\dot{R}_{i}(t)=\hspace{-0.1cm}\left\{
\begin{array}
[c]{ll}%
0 & \text{if }R_{i}(t)=0,\hspace{-0.2cm}\text{ }A_{i}\leq B_{i}P_{i}\left(
\mathbf{s}(t)\right) \\
A_{i}-B_{i}P_{i}\left(  \mathbf{s}(t)\right)  & \text{otherwise}%
\end{array}
\right.  \label{eq:multiDynR}%
\end{equation}
where we assume that initial conditions $R_{i}(0)$, $i=1,\ldots,M$, are given
and that $B_{i}>A_{i}>0$ (thus, the uncertainty strictly decreases when there
is perfect sensing $P_{i}\left(  \mathbf{s}(t)\right)  =1$).

\begin{figure}[h]
\begin{center}
\includegraphics[width=0.5\textwidth]{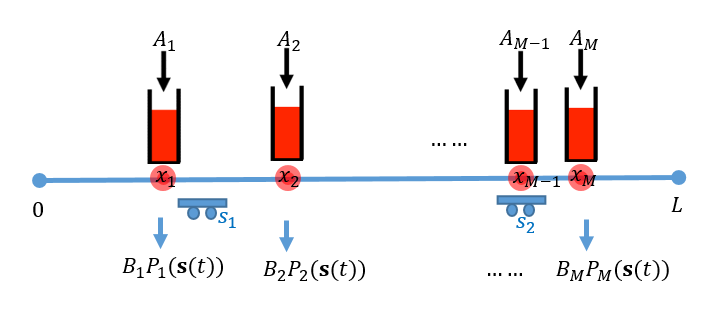}
\end{center}
\caption{\small A 1D polling model interpretation of problem \textbf{P1}}%
\label{fig:1DimModel}%
\end{figure}

Our goal is to control the movement of the $N$ agents through $u_{j}\left(
t\right)  $ in (\ref{eq:multiDynOfS}) so that the cumulative average uncertainty over
all targets $i=1,\ldots,M$ is minimized over a fixed time horizon $T$. Thus,
setting $\mathbf{u}\left(  t\right)  =\left[  u_{1}\left(  t\right)
,\ldots,u_{N}\left(  t\right)  \right]  $ we aim to solve the following
optimal control problem\textbf{ P1}:
\begin{equation}
\min_{\mathbf{u}\left(  t\right)  }\text{ \ }J=\frac{1}{T}\int_{0}^{T}%
\sum_{i=1}^{M}R_{i}(t)dt\label{eq:costfunction}%
\end{equation}
subject to the agent dynamics (\ref{eq:multiDynOfS}), uncertainty dynamics
(\ref{eq:multiDynR}), control constraint $|u_{j}(t)|\leq1$, $t\in\lbrack0,T]$,
and state constraints (\ref{eq:multiNoCross}). Figure \ref{fig:1DimModel} is a
polling model version for problem \textbf{P1} where each target is associated
with a \textquotedblleft virtual queue\textquotedblright\ where uncertainty
accumulates with inflow rate $A_{i}$. The service rate of this queue is
time-varying and given by $B_{i}P_{i}\left(  \mathbf{s}(t)\right)  $,
controllable through the agent position at time $t$. This interpretation is
convenient for characterizing the \emph{stability} of such a system over a
mission time $T$: For each queue, we may require that $\int_{0}^{T}A_{i}%
<\int_{0}^{T}B_{i}P_{i}(\mathbf{s}(t))dt$. Alternatively, we may require that
each queue becomes empty at least once over $[0,T]$. Note that this analogy
readily extends to two or three-dimensional settings. 

\section{Optimal control solution}
\label{sec:optimal}

In this section, we derive properties of the optimal control solution of
problem\textbf{ P1} and show that it can be reduced to a parametric
optimization problem. This will allow us to utilize an Infinitesimal
Perturbation Analysis (IPA) gradient estimation approach
\cite{cassandras2010perturbation} to find a complete optimal solution through
a gradient-based algorithm. We begin by defining the state vector
$\mathbf{x}(t)=[R_{1}(t),...R_{M}(t),s_{1}(t)...s_{N}(t)]$ and associated
costate vector $\lambda=[\lambda_{1}(t),...,\lambda_{M}(t),\lambda_{s_{1}%
}(t),...,\lambda_{s_{N}}(t)]$. As in \cite{cassandras2013optimal}, since the discontinuity in the dynamics of $R_{i}(t)$ in
(\ref{eq:multiDynR}), the optimal state trajectory may contain a boundary arc
when $R_{i}(t)=0$ for some $i$; otherwise, the state evolves in an interior
arc. Thus, we first analyze such an interior arc. Using (\ref{eq:multiDynOfS})
and (\ref{eq:multiDynR}), the Hamiltonian is
\begin{equation}
H(\mathbf{x},\lambda,\mathbf{u})=\sum_{i=1}^{M}R_{i}(t)+\sum_{i=1}^{M}\lambda
_{i}(t)\dot{R}_{i}(t)+\sum_{j=1}^{N}\lambda_{s_{j}}(t)u_{j}%
(t)\label{eq:Hamiltonian}%
\end{equation}
The costate dynamics are
\begin{equation}
\dot{\lambda_{i}}(t)=-\frac{\partial H}{\partial R_{i}(t)}=-1,\quad\lambda
_{i}(T)=0\label{eq:CostateDynTarget}%
\end{equation}%
\begin{equation}
\dot{\lambda}_{s_{j}}\left(  t\right)  =-\frac{\partial H}{\partial
s_{j}\left(  t\right)  }=\sum_{i=1}^{M}\lambda_{i}\left(  t\right)  B_{i}%
\frac{\partial P_{i}(\mathbf{s}(t))}{\partial s_{j}(t)},\quad\lambda_{s_{j}%
}(T)=0\label{eq:CostateDynAgent}%
\end{equation}
Applying the Pontryagin Minimum Principle to (\ref{eq:Hamiltonian}) with $\mathbf{u}^{\star}(t)$,
$t\in\lbrack0,T)$, denoting an optimal control, a necessary condition for
optimality is
\begin{equation}
H\left(  \mathbf{x}^{\star},\mathbf{\lambda}^{\star},\mathbf{u}^{\star
}\right)  =\min_{u_{j}\in\lbrack-1,1],\text{ }j=1,\ldots,N}H\left(
\mathbf{x},\mathbf{\lambda},\mathbf{u}\right)  \label{eq:PMP}%
\end{equation}
from which it immediately follows that%
\begin{equation}
u_{j}^{\ast}(t)=%
\begin{cases}
1 & \quad\text{if }\lambda_{s_{j}}(t)<0\\
-1 & \quad\text{if }\lambda_{s_{j}}(t)>0
\end{cases}
\label{eq:NecessaryConditionPMP}%
\end{equation}
Note that there exists a possibility that $\lambda_{s_{j}}\left(  t\right)
=0$ over some finite singular intervals\ \cite{bryson1975applied}, in which
case $u_{j}^{\ast}(t)$ may take values in $\{$ $-1,0,1\}$.

Similar to the case of the persistent monitoring problem studied in
\cite{cassandras2013optimal}, the complete solution requires solving the
costate equations (\ref{eq:CostateDynAgent}), which in turn involves the
determination of all points where $R_{i}(t)=0$, $i=1,\ldots,M$. This generally
involves the solution of a two-point-boundary-value problem. However, we will
next prove some structural properties of an optimal trajectory, based on which
we show that it is fully characterized by a set of parameters, thus reducing
the optimal control problem to a much simpler parametric optimization problem.

We begin by assuming that targets are ordered according to their location so
that $x_{1}<\cdots<x_{M}$. Let $r=\max_{j=1,\ldots,N}\{r_{j}\}$ and
$a=\max\{0,x_{1}-r\}$, $b=\min\{L,x_{M}+r\}$. Thus, if $s_{j}(t)<x_{1}-r$ or
$s_{j}(t)>x_{M}+r$, then it follows from (\ref{eq:SensingCap_1A}) that
$p_{ij}(s_{j}(t))=0$ for all targets $i=1,\ldots,M$. Clearly, this implies
that the effective mission space is $[a,b]$, i.e.,%
\begin{equation}
a\leq s_{j}(t)\leq b,\text{ \ \ }j=1,\ldots,N\label{eq:TrajConstrain}%
\end{equation}
imposing an additional state constraint for \textbf{P1}. We will show next
that on an optimal trajectory every agent is constrained to move within the
interval $[x_{1},x_{M}]$. This implies that every agent must switch its
direction no later than reaching the first or last target (possibly after
dwelling at the switching point for a finite time interval). To establish this
and subsequent results, we will make a technical assumption that no two events
altering the dynamics in this system can occur at the exact same time.

\textbf{Assumption 1}: \label{NoConcurrentEvent} Suppose that an agent
switches direction at $\theta\in\lbrack a,b]$. For any $j=1,\ldots,N,$
$i=1,\ldots,M,$ $t\in(0,T)$, and any $\epsilon>0$, if $s_{j}(t)=\theta$,
$s_{j}(t-\epsilon)>\theta$ or if $s_{j}(t)=\theta$, $s_{j}(t-\epsilon)<\theta
$, then either $R_{i}(\tau)>0$ for all $\tau\in\lbrack t-\epsilon,t]$ or
$R_{i}(\tau)=0$ for all $\tau\in\lbrack t-\epsilon,t]$.

\textbf{Proposition 1}: In an optimal trajectory, $x_{1}\leq s_{j}^{\ast
}(t)\leq x_{M}$, $t\in\lbrack0,T]$, $j=1,\ldots,N$.

\textbf{Proof}. We first prove that $s_{j}^{\ast}(t)\geq x_{1}$ for any agent
$j$. Suppose that $s_{j}^{\ast}(t_{0})=x_{1}$ and $u_{j}^{\ast}(t_{0})=-1$. In
view of (\ref{eq:TrajConstrain}), assume that agent $j$ reaches a point
$\theta\in\lbrack a,x_{1})$ at time $t_{1}>t_{0}$ where it switches direction;
we will show that $\theta\notin\lbrack a,x_{1})$ using a contradiction
argument. There are two cases to consider:

\emph{Case 1}: $\theta=a$. Assuming $s_{j}^{\ast}(t_{1})=a$, we first show
that $\lambda_{s_{j}}^{\ast}(t_{1}^{-})=0$ by a contradiction argument. If
$\lambda_{s_{j}}^{\ast}(t_{1}^{-})\neq0$, recall that $u_{j}^{\ast}(t_{1}%
^{-})=-1$, therefore $\lambda_{s_{j}}^{\ast}(t_{1}^{-})>0$ from
(\ref{eq:NecessaryConditionPMP}). Since the constraint $a-s_{j}(t)\leq0$ is
active, $\lambda_{s_{j}}^{\ast}(t)$ may experience a discontinuity so that
\begin{equation}
\lambda_{s_{j}}^{\ast}(t_{1}^{-})=\lambda_{s_{j}}^{\ast}(t_{1}^{+})-\pi_{j}%
\end{equation}
where $\pi_{j}\geq0$ is a scalar multiplier associated with the constraint
$a-s_{j}(t)\leq0$. It follows that $\lambda_{s_{j}}^{\ast}(t_{1}^{+}%
)=\lambda_{s_{j}}^{\ast}(t_{1}^{-})+\pi_{j}>0$. Since the Hamiltonian in
(\ref{eq:Hamiltonian}) and the constraint $a-s_{j}(t)\leq0$ are not explicit
functions of time, we have \cite{bryson1975applied} $H^{\ast}(\mathbf{x}%
(t_{1}^{-}),\mathbf{\lambda}(t_{1}^{-}),\mathbf{u}(t_{1}^{-}))=H^{\ast
}(\mathbf{x}(t_{1}^{+}),\lambda(t_{1}^{+}),\mathbf{u}(t_{1}^{+}))$ which,
under Assumption 1, reduces to
\begin{equation}
\lambda_{s_{j}}^{\ast}(t_{1}^{-})u_{j}^{\ast}(t_{1}^{-})=\lambda_{s_{j}}%
^{\ast}(t_{1}^{+})u_{j}^{\ast}(t_{1}^{+})\label{eq:lemma1_1}%
\end{equation}
Recall that $\lambda_{s_{j}}^{\ast}(t_{1}^{-})u_{j}^{\ast}(t_{1}^{-})<0$.
However, $u_{j}^{\ast}(t_{1}^{+})\geq0$ (since the agent switches control),
therefore $\lambda_{s_{j}}^{\ast}(t_{1}^{+})u_{j}^{\ast}(t_{1}^{+})\geq0$
which violates (\ref{eq:lemma1_1}). This contradiction implies that
$\lambda_{s_{j}}^{\ast}(t_{1}^{-})=0$. Recalling (\ref{eq:SensingCap_NA}) and
(\ref{eq:CostateDynAgent}), we get $\dot{\lambda}_{s_{j}}^{\ast}(t_{1}^{-}%
)=\sum_{i=1,R_{i}\neq0}^{M}\lambda_{i}^{\ast}(t_{1}^{-})\frac{B_{i}}{r_{j}}%
\prod_{d\neq j}[1-p_{id}(s_{d}^{\ast}(t_{1}^{-}))]$. Under Assumption 1, there exists
$\delta>0$ such that during interval $(t_{1}-\delta,t_{1})$, no $R_{i}%
(t)\geq0$ becomes active, hence no $\lambda_{i}^{\ast}(t)$ encounters a jump for
$i=1,\ldots,M$ and it follows from (\ref{eq:CostateDynTarget}) that $\lambda
_{i}^{\ast}(t)>0$. Moreover, $p_{id}(s_{d}^{\ast}(t))\neq1$ for at least some $d\neq j$
since we have assumed that $M>N$. Thus, we have $\dot{\lambda}_{s_{j}}^{\ast}(t)>0,$
for all $t\in(t_{1}-\delta,t_{1})$. However, since agent $j$ is approaching
$a$, there exists some $\delta^{\prime}<\delta$, such that $u_{j}^{\ast}(t)=-1$ for
all $t\in(t_{1}-\delta^{\prime},t_{1})$, and $\lambda_{s_{j}}^{\ast}(t)\geq0$. Thus
for $t\in(t_{1}-\delta^{\prime},t_{1})$, we have $\lambda_{s_{j}}^{\ast}(t)\geq0$ and
$\dot{\lambda}_{s_{j}}^{\ast}(t)>0$. This contradicts the established fact that
$\lambda_{s_{j}}^{\ast}(t_{1}^{-})=0$. We conclude that $\theta\neq a$.

\emph{Case 2}: $\theta\in(a,x_{1})$. Assuming $s_{j}^{\ast}(t_{1})=\theta$, we
still have $u_{j}^{\ast}(t_{1}^{-})=-1$, $u_{j}^{\ast}(t_{1}^{+})\geq0$. Since
the Hamiltonian (\ref{eq:Hamiltonian}) is not an explicit function of time, we
have $H^{\ast}(\mathbf{x}(t_{1}^{-}),\mathbf{\lambda}(t_{1}^{-}),\mathbf{u}%
(t_{1}^{-}))=H^{\ast}(\mathbf{x}(t_{1}^{+}),\lambda(t_{1}^{+}),\mathbf{u}%
(t_{1}^{+}))$ which leads to (\ref{eq:lemma1_1}) under Assumption 1. First, we
assume $\lambda_{s_{j}}^{\ast}(t_{1}^{-})\neq0$. Since $u_{j}^{\ast}(t_{1}%
^{-})<0$, we have $\lambda_{s_{j}}^{\ast}(t_{1}^{-})>0$ and the left hand side
of (\ref{eq:lemma1_1}) gives $\lambda_{s_{j}}^{\ast}(t_{1}^{-})u_{j}^{\ast}(t_{1}%
^{-})<0$. On the other hand, in order to satisfy (\ref{eq:lemma1_1}), we must
have $u_{j}^{\ast}(t_{1}^{+})>0$ and $\lambda_{s_{j}}^{\ast}(t_{1}^{+})<0$. However,
if $\lambda_{s_{j}}^{\ast}(t_{1}^{-})>0$ and $\lambda_{s_{j}}^{\ast}(t_{1}^{+})<0$,
then either $\dot{\lambda}_{s_{j}}^{\ast}(t_{1})<0$ and $\lambda_{s_{j}}^{\ast}(t_{1})=0$,
or $\lambda_{s_{j}}^{\ast}(t)$ experiences a discontinuity at $t_{1}$. We show that
neither condition is feasible. The first one violates our assumption that
$\lambda_{s_{j}}^{\ast}(t_{1})\neq0$, while the second one is not feasible since at
$t=t_{1}$ the constraint $a-s_{j}(t)\leq0$ is not active. This implies that
$\lambda_{s_{j}}^{\ast}(t_{1}^{-})=0$. Again, under Assumption 1, the same
argument as in \emph{Case 1} can be used to show that $\lambda_{s_{j}}%
^{\ast}(t)\geq0$ and $\dot{\lambda}_{s_{j}}^{\ast}(t)>0$ for all $t\in(t_{1}-\delta^{\prime
},t_{1})$. This contradicts the established fact that $\lambda_{s_{j}}%
^{\ast}(t_{1}^{-})=0$ and we conclude that $\theta\notin(a,x_{1})$.

Combining both cases, we conclude that $\theta\notin\lbrack0,x_{1})$, which
implies that $s_{j}^{\ast}(t)\geq x_{1}$. The same line of argument can be
used to show that $s_{j}^{\ast}(t)\leq x_{M}$.$\blacksquare$

Proposition 1, in conjunction with (\ref{eq:NecessaryConditionPMP}), leads to the
conclusion that the optimal control consists of each agent moving with maximal
speed in one direction until it reaches a point in the interval $[x_{1}%
,x_{M}]$ where it switches direction. However, the exclusion of the case
$\lambda_{s_{j}}(t)=0$ allows the possibility of \emph{singular arcs} along
the optimal trajectory, defined as intervals $[t_{1},t_{2}]$ such that
$\lambda_{s_{j}}(t)=0$ for all $t\in\lbrack t_{1},t_{2}]$ and $\lambda_{s_{j}%
}(t_{1}^{-})\neq0$, $\lambda_{s_{j}}(t_{2}^{+})\neq0$. The next result
establishes the fact that we can exclude singular arcs from an agent's
trajectory while this agent has no target in its sensing range.

\textbf{Lemma 1}: If $|s_{j}(t)-x_{i}|>r_{j}$ for any $i=1,\ldots,M$, then
$u_{j}^{\ast}(t)\neq0$.

\textbf{Proof}. We proceed with a contradiction argument. Suppose that
$u_{j}^{\ast}(t)=0$ for $t\in\lbrack t_{1},t_{2}]$ such that $|s_{j}^{\ast
}(t_{1})-x_{i}|>r_{j}$ for all $i=1,\ldots,M$ and that $u_{j}^{\ast}(t)\neq0$
(without loss of generality, let $u_{j}^{\ast}(t)=1$) for $t>t_{2}$ so that
$|s_{j}^{\ast}(t_{3})-x_{i}|=r_{j}$ for some $i=1,\ldots,M$ and $|s_{j}^{\ast
}(t_{3}+\Delta)-x_{i}|<r_{j}$ for $t_{3}+\Delta>t_{3}>t_{2}$. In other words,
agent $j$ eventually reaches a target $i$ that it can sense at $t=t_{3}$.
Assume that $u_{j}^{\ast}(t)$, $t\in\lbrack t_{1},t_{3}+\Delta]$ is replaced
by $u_{j}^{\prime}(t)$ as follows: $u_{j}^{\prime}(t)=1$ for $t\in\lbrack
t_{1},t_{3}+\Delta+t_{1}-t_{2}]$ and $u_{j}^{\prime}(t)=0$ for $t\in
(t_{3}+\Delta+t_{1}-t_{2},t_{3}+\Delta]$. In other words, the agent moves to
reach $s_{j}^{\prime}(t_{3}+\Delta+t_{1}-t_{2})=s_{j}^{\ast}(t_{3}+\Delta)$
and then stops. The two controls are thereafter identical. Then, referring to
(\ref{eq:costfunction}) we have $\int_{t_{3}+\Delta+t_{1}-t_{2}}^{t_{3}%
+\Delta}R_{i}^{\prime}(t)dt\leq\int_{t_{3}+\Delta+t_{1}-t_{2}}^{t_{3}+\Delta
}R_{i}^{\ast}(t)dt$ since under $u_{j}^{\prime}(t)$ the agent may decrease
$R_{i}(t)$ over $[t_{3}+\Delta+t_{1}-t_{2},t_{3}]$ whereas under $u_{j}^{\ast
}(t)$ this is impossible since $|s_{j}^{\ast}(t)-x_{i}|>r_{j}$ over this time
interval. Since the cost in (\ref{eq:costfunction}) is the same over
$[0,t_{3}+\Delta+t_{1}-t_{2})$ and $(t_{3}+\Delta,T],$ it follows that
$u_{j}^{\ast}(t)=0$ when $|s_{j}(t)-x_{i}|>r_{j}$ cannot be optimal unless
$u_{j}^{\ast}(t)=0$ for all $t\in\lbrack0,T]$, i.e., the agent never moves and
never senses any target, in which case the cost under $u_{j}^{^{\prime}}(t)$
is still no higher than that under $u_{j}^{\ast}(t)$.$\blacksquare$

Based on Lemma 1, we conclude that singular arcs in an agent's trajectory may
occur only while it is sensing a target. Intuitively, this indicates that it
may be optimal for an agent to stop moving and dwell in the vicinity of one or
more targets that it can sense so as to decrease the associated uncertainty
functions to an adequate level before it proceeds along the mission space. The
next lemma establishes the fact that if the agent is visiting an
\emph{isolated} target and experiences a singular arc, then the corresponding
optimal control is $u_{j}^{\ast}(t)=0$. An isolated target with position
$x_{i}$ is defined to be one that satisfies $|x_{i}-x_{j}|>2r,$for all $j\neq
i$ where $r$ was defined earlier as $r=\max_{j=1,\ldots,N}\{r_{j}\}$.
Accordingly, the subset $I\sqsubseteq\{1,\ldots,M\}$ of isolated targets is
defined as%
\begin{equation}
I=\{i:|x_{i}-x_{j}|>2r,j\neq i\in\{1,\makebox[1em][c]{.\hfil.\hfil.}%
,M\},r=\hspace{-0.3cm}\max_{j=1,\makebox[1em][c]{.\hfil.\hfil.},N}\{r_{j}\}\}
\label{eq:isolated_target}%
\end{equation}

\textbf{Lemma 2}: Let $|s_{j}^{\ast}(t)-x_{k}|<r_{j}$ for some $j=1,\ldots,N$
and isolated target $k\in I$. If $\lambda_{s_{j}}^{\ast}(t)=0$, $t\in\lbrack
t_{1},t_{2}]$, then $u_{j}^{\ast}(t)=0$.

\textbf{Proof}. The proof is along the same line as Proposition III.3 in
\cite{cassandras2013optimal}. Assume that $\lambda_{s_{j}}^{\ast}(t)=0$ over a
singular arc $[t_{1},t_{2}]$. Let $H^{\ast}\equiv$ $H(\mathbf{x}^{\ast
}\mathbf{,\lambda}^{\ast}\mathbf{,u}^{\ast})$. Since the Hamiltonian along an
optimal trajectory is a constant, we have $\frac{dH^{\ast}}{dt}=0$. Therefore,
recalling (\ref{eq:Hamiltonian}),
\begin{gather*}
\frac{dH^{\ast}}{dt}=\sum_{i=1}^{M}\Big[\dot{R}_{i}^{\ast}(t)+\dot{\lambda
}_{i}^{\ast}(t)\dot{R}_{i}^{\ast}(t)+\lambda_{i}^{\ast}(t)\ddot{R}_{i}^{\ast
}(t)\Big]\\
+\sum_{j=1}^{N}\Big[\dot{\lambda}_{s_{j}}^{\ast}(t)u_{j}^{\ast}(t)+\lambda
_{s_{j}}^{\ast}(t)\dot{u}_{j}^{\ast}(t)\Big]=0
\end{gather*}
and since $\dot{\lambda_{i}^{\ast}}(t)=-1$ from (\ref{eq:CostateDynTarget}), this reduces to
\begin{equation}
\frac{dH^{\ast}}{dt}=\sum_{i=1}^{M}\lambda_{i}^{\ast}(t)\ddot{R}_{i}^{\ast
}(t)+\sum_{j=1}^{N}\Big[\dot{\lambda}_{s_{j}}^{\ast}(t)u_{j}^{\ast}%
(t)+\lambda_{s_{j}}^{\ast}(t)\dot{u}_{j}^{\ast}(t)\Big]=0 \label{eq:Lemma3_1}%
\end{equation}
Define $S(t)=\{j|\lambda_{s_{j}}(t)=0,\dot{\lambda}_{s_{j}}(t)=0\}$ as the set
of agents in singular arcs at $t$ and $\bar{S}(t)$ as the set of all remaining
agents. If $j\in S(t)$, then $\dot{\lambda}_{s_{j}}^{\ast}(t)u_{j}^{\ast
}(t)+\lambda_{s_{j}}^{\ast}(t)\dot{u}_{j}^{\ast}(t)=0$. If $j\in\bar{S}(t)$,
then $\lambda_{s_{j}}^{\ast}(t)\dot{u}_{j}^{\ast}(t)=0$ since $u_{j}^{\ast
}(t)=\pm1$ and $\dot{u}_{j}^{\ast}(t)=0$. Therefore, we rewrite
(\ref{eq:Lemma3_1}) as
\begin{equation}
\frac{dH^{\ast}}{dt}=\sum_{i=1}^{M}\lambda_{i}^{\ast}(t)\ddot{R}_{i}^{\ast
}(t)+\sum_{j\in\bar{S}(t)}\dot{\lambda}_{s_{j}}^{\ast}(t)u_{j}^{\ast}(t)=0
\label{eq:Lemma3_2}%
\end{equation}
Recalling (\ref{eq:multiDynR}), when $R_{i}(t)\neq0$, we have $\dot{R}%
_{i}=A_{i}-B_{i}\Big(1-\prod_{n=1}^{N}\Big[1-p_{ij}\big(s_{j}%
(t)\big)\Big]\Big)$. Therefore,
\begin{align}
\ddot{R}_{i}^{\ast}(t)=  &  \frac{d}{dt}\dot{R}_{i}^{\ast}(t)\nonumber\\
&  =-\sum_{j=1}^{N}u_{j}^{\ast}(t)B_{i}\frac{\partial p_{ij}(s_{j}^{\ast}%
(t))}{\partial s_{j}^{\ast}}\prod_{d\neq j}\Big[1-p_{id}(s_{d}^{\ast}(t))\Big]
\end{align}

Moreover, from (\ref{eq:CostateDynAgent}), we have
\begin{equation}
\dot{\lambda}_{s_{j}}^{\ast}(t)=\sum_{i=1,R_{i}\neq0}^{M}\lambda_{i}^{\ast
}(t)B_{i}\frac{\partial p_{ij}(s_{j}^{\ast}(t))}{\partial s_{j}^{\ast}}%
\prod_{d\neq j}\Big[1-p_{id}(s_{d}^{\ast}(t))\Big]\label{eq:Lemma3_3}%
\end{equation}
Combining \eqref{eq:Lemma3_2}-\eqref{eq:Lemma3_3}, we get
\begin{align}
&  \frac{dH^{\ast}}{dt}= \nonumber\\
&  -{\sum\limits_{\substack{{i=1} \\{R_{i}\neq0}}}^{M}}\sum_{j=1}^{N}u_{j}^{\ast}(t)\lambda_{i}^{\ast
}(t)B_{i}\frac{\partial p_{ij}(s_{j}^{\ast}(t))}{\partial s_{j}^{\ast}}%
\prod_{d\neq j}\Big[1-p_{id}(s_{d}^{\ast}(t))\Big]\nonumber\\
&  +\sum_{j\in\bar{S}(t)}{\sum\limits_{\substack{{i=1} \\{R_{i}\neq0}}}^{M}}u_{j}^{\ast}(t)\lambda
_{i}^{\ast}(t)B_{i}\frac{\partial p_{ij}(s_{j}^{\ast}(t))}{\partial
s_{j}^{\ast}}\prod_{d\neq j}\Big[1-p_{id}(s_{d}^{\ast}(t))\Big]\nonumber\\
&  =\sum_{j\in S(t)}{\sum\limits_{\substack{{i=1} \\{R_{i}\neq0}}}^{M}} u_{j}^{\ast}(t)\lambda_{i}^{\ast
}(t)B_{i}\frac{\partial p_{ij}(s_{j}^{\ast}(t))}{\partial s_{j}^{\ast}}%
\prod_{d\neq j}\Big[1-p_{id}(s_{d}^{\ast}(t))\Big]\nonumber\\
&  =0\label{eq:Lemma3_4}%
\end{align}
Since we have assumed that $|s_{j}^{\ast}(t)-x_{k}|<r_{j}$ and $k$ is an
isolated target, it follows that $p_{kj}(s_{j}^{\ast}(t))\neq0$ and
$p_{ij}(s_{j}(t))=0$ if $i\neq k$. Therefore, $\frac{\partial p_{kj}%
(s_{j}^{\ast}(t))}{\partial s_{j}^{\ast}}\neq0$ and $\frac{\partial
p_{ij}(s_{j}^{\ast}(t))}{\partial s_{j}^{\ast}}=0$ for all $i\neq k$ and
(\ref{eq:Lemma3_4}) reduces to
\begin{equation}
\sum_{j\in S(t)}u_{j}^{\ast}(t)\lambda_{k}^{\ast}(t)B_{i}\frac{\partial
p_{kj}(s_{j}^{\ast}(t))}{\partial s_{j}^{\ast}}\prod_{d\neq j}\Big[1-p_{kd}%
(s_{d}^{\ast}(t))\Big]=0\label{eq:Lemma3_5}%
\end{equation}
Observe that, from (\ref{eq:CostateDynTarget}), $\lambda_{i}(t)>0$ when
$R_{i}(t)\neq0,t<T$. In addition $B_{i}>0$ and $\prod_{d\neq j}[1-p_{kd}%
(s_{d}^{\ast}(t))]\neq0$. Therefore, to satisfy (\ref{eq:Lemma3_5}) for all
$t\in\lbrack t_{1},t_{2}]$, we must have $u_{j}^{\ast}(t)=0,$ for all $j\in
S(t)$.$\blacksquare$

We can further establish the fact that if an agent $j$ experiences a singular
arc while sensing an isolated target $k$, then the optimal point to stop is
such that $s_{j}^{\ast}(t)=x_{k}$.

\textbf{Proposition 2}: Let $|s_{j}^{\ast}(t)-x_{k}|<r_{j}$ for some
$j=1,\ldots,N$ and isolated target $k\in I$. If $\lambda_{s_{j}}^{\ast}(t)=0$,
$t\in\lbrack t_{1},t_{2}]$, and $u_{j}^{\ast}(t_{1}^{-})=u_{j}^{\ast}%
(t_{2}^{+})$, then $s_{j}^{\ast}(t)=x_{k}$, $t\in\lbrack t_{1},t_{2}]$.

\textbf{Proof}. By Lemma 2, we know that $u_{j}^{\ast}(t)=0$, $t\in\lbrack
t_{1},t_{2}]$. We use a contradiction argument similar to the one used in
Lemma 1 to show that $s_{j}^{\ast}(t)=x_{k}$, $t\in\lbrack t_{1},t_{2}]$.
Suppose that $u_{j}^{\ast}(t_{1}^{-})=1$ (without loss of generality) and that
$s_{j}^{\ast}(t)=x_{k}-\Delta<x_{k}$. Note that at the end of the singular arc
$u_{j}^{\ast}(t_{2}^{+})=1$ since $u_{j}^{\ast}(t_{1}^{-})=u_{j}^{\ast}%
(t_{2}^{+})$. This implies that $s_{j}^{\ast}(t_{2}+\Delta)=x_{k}.$ Assume
that $u_{j}^{\ast}(t)$, $t\in\lbrack t_{1},t_{2}+\Delta]$ is replaced by
$u_{j}^{\prime}(t)$ as follows: $u_{j}^{\prime}(t)=1$ for $t\in\lbrack
t_{1},t_{1}+\Delta]$ and $u_{j}^{\prime}(t)=0$ for $t\in(t_{1}+\Delta
,t_{2}+\Delta]$. In other words, the agent moves to reach $s_{j}^{\prime
}(t_{1}+\Delta)=s_{j}^{\ast}(t_{2}+\Delta)=x_{k}$ and then stops. The two
controls are thereafter identical. Then, referring to (\ref{eq:costfunction})
we have $\int_{t_{1}}^{t_{2}+\Delta}R_{i}^{\prime}(t)dt<\int_{t_{1}}%
^{t_{2}+\Delta}R_{i}^{\ast}(t)dt$ since $\dot{R}_{i}^{\ast}(t)<\dot{R}%
_{i}^{\prime}(t)$ due to (\ref{eq:multiDynR}) and the fact that $p_{kj}%
(s_{j}(t))$ is monotonically decreasing in $|s_{j}(t)-x_{k}|$. Since the cost
in (\ref{eq:costfunction}) is the same over $[0,t_{1})$ and $(t_{2}%
+\Delta,T],$ it follows that $s_{j}^{\ast}(t)=x_{k}-\Delta$ cannot be optimal.
The same argument holds for any $\Delta>0$, leading to the conclusion that
$s_{j}^{\ast}(t)=x_{k}$, $t\in\lbrack t_{1},t_{2}]$. A similar argument also
applies to the case $s_{j}^{\ast}(t)=x_{k}+\Delta>x_{k}$.$\blacksquare$

Finally, we consider the case with the state constraint
(\ref{eq:multiNoCross}). We can then prove that this constraint is never
active on an optimal trajectory, i.e., agents reverse their directions before
making contact with any other agent.

\textbf{Proposition 3:} Under the constraint $s_{j}(t)\leq s_{j+1}(t)$, on an
optimal trajectory, $s_{j}(t)\neq s_{j+1}(t)$ for all $t\in(0,T),$ $j=1...N$.

\textbf{Proof}. The proof is almost identical to that of Proposition III.4 in
\cite{cassandras2013optimal} and is, therefore, omitted.$\blacksquare$

The above analysis, including Propositions 1-3, fully characterize the
structure of the optimal control as consisting of intervals in $[0,T]$ where
$u_{j}^{\ast}(t)\in\{-1,0,1\}$ depending entirely on the sign of
$\lambda_{s_{j}}(t)$. Based on this analysis, we can parameterize \textbf{P1}
so that the cost in (\ref{eq:costfunction}) depends on a set of $(i)$\textit{
}\emph{switching points} where an agent switches its control from
$u_{j}(t)=\pm1$ to $\mp1$ or possibly $0$, and $(ii)$\textit{ }\emph{dwelling
times} if an agent switches from $u_{j}(t)=\pm1$ to $0$. In other words, the
optimal trajectory of each agent $j$ is totally characterized by two parameter
vectors: switching points $\bm\theta_{j}=[\theta_{j1},\theta_{j2}%
...\theta_{j\Gamma}]$ and dwelling times $\bm\omega_{j}=[\omega_{j1}%
,\omega_{j2}...\omega_{j\Gamma^{\prime}}]$ where $\Gamma$ and $\Gamma^{\prime
}$ are prior parameters depending on the given time horizon. This defines a
hybrid system with state dynamics (\ref{eq:multiDynOfS}), (\ref{eq:multiDynR}%
). The dynamics remain unchanged in between events that cause them to change,
i.e., the points $\theta_{j1},\ldots,\theta_{j\Gamma}$ above and instants when
$R_{i}(t)$ switches from $>0$ to $0$ or vice versa. Therefore, the overall
cost function (\ref{eq:costfunction}) can be parametrically expressed as
$J(\bm\theta,\bm\omega)$ and rewritten as the sum of costs over corresponding
interevent intervals over a given time horizon:%
\begin{equation}
J(\bm\theta,\bm\omega)=\frac{1}{T}\sum_{k=0}^{K}\int_{\tau_{k}(\bm\theta
,\bm\omega)}^{\tau_{k+1}(\bm\theta,\bm\omega)}\sum_{i=1}^{M}R_{i}%
(t)dt\label{eq:ParametricObj}%
\end{equation}
where $\tau_{k}$ is the $k$-th event time. This will allow us to apply IPA to
determine a gradient $\nabla J(\bm\theta,\bm\omega)$ with respect to these
parameters and apply any standard gradient-based optimization algorithm to
obtain a (locally) optimal solution.

\section{Infinitesimal Perturbation Analysis}
\label{sec:IPA}
As concluded in the previous section, optimal agent trajectories may be
selected from the family $\{{\bm s(\bm\theta},\bm\omega,t,{\bm s}_{0})\}$ with
parameter vectors $\bm\theta$ and $\bm\omega$ and a given initial condition
${\bm s}_{0}$. Along these trajectories, agents are subject to dynamics
(\ref{eq:multiDynOfS}) and targets are subject to (\ref{eq:multiDynR}). An
\emph{event} (e.g., an agent stopping at some target $x_{i}$) occurring at
time $\tau_{k}({\bm\theta},\bm\omega)$ triggers a switch in these state
dynamics. IPA specifies how changes in $\bm\theta$ and $\bm\omega$ influence
the state ${\bm s(\bm\theta},\bm\omega,t,{\bm s}_{0})$, as well as event times
$\tau_{k}({\bm\theta},\bm\omega)$, $k=1,2,\ldots$, and, ultimately the cost
function (\ref{eq:ParametricObj}). We briefly review next the IPA framework for
general stochastic hybrid systems as presented in
\cite{cassandras2010perturbation}.

Let $\{\tau_{k}(\theta)\}$, $k=1,\ldots,K$, denote the occurrence times of all
events in the state trajectory of a hybrid system with dynamics $\dot
{x}\ =\ f_{k}(x,\theta,t)$ over an interval $[\tau_{k}(\theta),\tau
_{k+1}(\theta))$, where $\theta\in\Theta$ is some parameter vector and
$\Theta$ is a given compact, convex set. For convenience, we set $\tau_{0}=0$
and $\tau_{K+1}=T$. We use the Jacobian matrix notation: $x^{\prime}%
(t)\equiv\frac{\partial x(\theta,t)}{\partial \theta}$ and $\tau
_{k}^{\prime}\equiv\frac{\partial\tau_{k}(\theta)}{\partial \theta}$, for
all state and event time derivatives. It is shown in
\cite{cassandras2010perturbation} that
\begin{equation}
\frac{d}{dt}x^{\prime}(t)=\frac{\partial f_{k}(t)}{\partial x}x^{\prime
}(t)+\frac{\partial f_{k}(t)}{\partial\theta}\label{eq:IPA_1}%
\end{equation}
for $t\in\lbrack\tau_{k},\tau_{k+1})$ with boundary condition:
\begin{equation}
x^{\prime}(\tau_{k}^{+})=x^{\prime}(\tau_{k}^{-})+[f_{k-1}(\tau_{k}^{-}%
)-f_{k}(\tau_{k}^{+})]\tau_{k}^{\prime}\label{eq:IPA_2}%
\end{equation}
for $k=0,...K$. In order to complete the evaluation of $x^{\prime}(\tau
_{k}^{+})$ in (\ref{eq:IPA_2}), we need to determine $\tau_{k}^{\prime}$. If the
event at $\tau_{k}$ is \emph{exogenous}, $\tau_{k}^{\prime}=0$. However, if
the event is \emph{endogenous}, there exists a continuously differentiable
function $g_{k}:\mathbb{R}^{n}\times\Theta\rightarrow\mathbb{R}$ such that
$\tau_{k}\ =\ \min\{t>\tau_{k-1}\ :\ g_{k}\left(  x\left(  \theta,t\right)
,\theta\right)  =0\}$ and
\begin{equation}
\tau_{k}^{\prime}=-[\frac{\partial g_{k}}{\partial x}f_{k}(\tau_{k}^{-}%
)]^{-1}(\frac{\partial g_{k}}{\partial \theta}+\frac{\partial g_{k}%
}{\partial x} x^{\prime}(\tau_{k}^{-}))\label{eq:IPA_3}%
\end{equation}
as long as $\frac{\partial g_{k}}{\partial x}f_{k}(\tau_{k}^{-})\neq0$.
(details may be found in \cite{cassandras2010perturbation}).

Denote the time-varying cost along a given trajectory as $L( x, \theta
,t)$, so the cost in the $k$-th interevent interval is $J_{k}(x, \theta
)=\int_{\tau_{k}}^{\tau_{k+1}}L( x,\theta,t )dt$ and the total cost is
$J(x, \theta)=\sum_{k=0}^{K}J_{k}( x,\theta)$. Differentiating and
applying the Leibnitz rule with the observation that all terms of the form
$L( x(\tau_{k}),\theta,\tau_{k})\tau_{k}^{\prime}$ are mutually canceled
with $\tau_{0}=0,\tau_{K+1}=T$ fixed, we obtain%
\begin{align}
\frac{\partial J(x,\theta)}{\partial \theta} &  =\sum_{k=0}^{K}%
\frac{\partial}{\partial\theta}\int_{\tau_{k}}^{\tau_{k+1}}L(
x,\theta,t)dt\nonumber\\
&  =\sum_{k=0}^{K}\int_{\tau_{k}}^{\tau_{k+1}}\frac{\partial L(
x,\theta,t)}{\partial x} x^{\prime}(t)+\frac{\partial L(
x,\theta,t)}{\partial\theta}dt\label{eq:GeneralGradientParametricObj}%
\end{align}
In our setting, we have $L( x,\theta,t)=\sum_{i=1}^{M}R_{i}(t)$ from
(\ref{eq:ParametricObj}), which is not an explicit function of the state
$\mathbf{x}(t)=[R_{1}(t),...R_{M}(t),s_{1}(t)...s_{N}(t)]$. Thus, the gradient
$\nabla J(\bm \theta,\bm \omega)=[\frac{\partial J(\bm\theta,\bm\omega
)}{\partial\bm\theta},\frac{\partial J(\bm\theta,\bm\omega)}%
{\partial\bm\omega}]^{\text{T}}$ reduces to
\begin{equation}
\nabla J(\bm\theta,\bm\omega)=\frac{1}{T}\sum_{k=0}^{K}\sum_{i=1}^{M}%
\int_{\tau_{k}(\bm\theta,\bm\omega)}^{\tau_{k+1}(\bm\theta,\bm\omega)}\nabla
R_{i}(t)dt\label{eq:GradientParametricObj}%
\end{equation}
where $\nabla R_{i}(t)=[\frac{\partial R_{i}(t)}{\partial\bm\theta}%
,\frac{\partial R_{i}(t)}{\partial\bm\omega}]^{\text{T}}$.

Applying (\ref{eq:IPA_1})(\ref{eq:IPA_2})(\ref{eq:IPA_3}), we can evaluate $\nabla
R_{i}(t)$. In contrast to \cite{cassandras2013optimal}, in our problem agents
are allowed to dwell on every target and IPA will optimize these dwelling
times. Therefore, we need to consider all possible forms of control sequences:
$(i)$ $\pm1\rightarrow0$, $(ii)$ $0\rightarrow\pm1$, and $(iii)$
$\pm1\rightarrow\mp1$. We can then obtain from (\ref{eq:multiDynR}):%

\begin{equation}
\frac{\partial R_{i}(t)}{\partial\theta_{j\xi}}=\frac{\partial R_{i}(\tau
_{k}^{+})}{\partial\theta_{j\xi}}-%
\begin{cases}
0\quad\text{if }R_{i}(t)=0,A_{i}<B_{i}P_{i}(\mathbf{s}(t))\\
G\frac{\partial s_{j}(\tau_{k}^{+})}{\partial\theta_{j\xi}}(t-\tau_{k}%
)\quad\text{otherwise}%
\end{cases}
\label{eq:pRptheta}%
\end{equation}

\begin{equation}
\frac{\partial R_{i}(t)}{\partial\omega_{j\xi}}=\frac{\partial R_{i}(\tau
_{k}^{+})}{\partial\omega_{j\xi}}-%
\begin{cases}
0\quad\text{if }R_{i}(t)=0,A_{i}<B_{i}P_{i}(\mathbf{s}(t))\\
G\frac{\partial s_{j}(\tau_{k}^{+})}{\partial\omega_{j\xi}}(t-\tau_{k}%
)\quad\text{otherwise}%
\end{cases}
\label{eq:pRpomega}%
\end{equation}
where $G=B_{i}\prod_{d\neq j}\big[1-p_{i}(s_{d}(t))\big]\frac{\partial
p_{i}(s_{j})}{\partial s_{j}}$ and $\frac{\partial p_{i}(s_{j})}{s_{j}}%
=\pm\frac{1}{r_{j}}$.

First, let us consider the events that cause switches in $\dot{R}_{i}(t)$ in
(\ref{eq:multiDynR}) at time $\tau_{k}$. For these events, the dynamics of
$s_{j}(t)$ are continuous so that $\nabla s_{j}(\tau_{k}^{-})=\nabla
s_{j}(\tau_{k}^{+})$. For target $i$,%

\begin{equation}
\nabla R_{i}(\tau_{k}^{+})=%
\begin{cases}
\nabla R_{i}(\tau_{k}^{-}) & \text{\textbf{if }}\dot{R}_{i}(\tau_{k}^{-})=0,\\
& \dot{R}_{i}(\tau_{k}^{+})=A_{i}-B_{i}P_{i}({\bm s}(\tau_{k}^{+})).\\
0 & \text{\textbf{if }}\dot{R}_{i}(\tau_{k}^{-})=A_{i}-B_{i}P_{i}({\bm s}%
(\tau_{k}^{-})),\\
& \dot{R}_{i}(\tau_{k}^{+})=0.
\end{cases}
\end{equation}

Second, let us consider events that cause switches in $\dot{s}_{j}%
(t)=u_{j}(t)$ at time $\tau_{k}$. For these events, the dynamics of $R_{i}(t)$
are continuous so that $\nabla R_{i}(\tau_{k}^{-})=\nabla R_{i}(\tau_{k}^{+}%
)$. In order to evaluate (\ref{eq:pRptheta}) and (\ref{eq:pRpomega}), we need
$\frac{\partial s_{j}(\tau_{k}^{+})}{\partial\theta_{j\xi}}$ and
$\frac{\partial s_{j}(\tau_{k}^{+})}{\partial\omega_{j\xi}}$. Clearly, these
are not affected by future events and we only have to consider the current and
prior control switches. Let $\theta_{j\xi}$ and $\omega_{j\xi}$ be the current
switching point and dwelling time. Again, applying
(\ref{eq:IPA_1})(\ref{eq:IPA_2})(\ref{eq:IPA_3}), we have

\emph{Case 1}: $u_{j}(\tau_{k}^{-})=\pm1,u_{j}(\tau_{k}^{+})=0$
\begin{align}
&  \frac{\partial s_{j}}{\partial\theta_{jl}}(\tau_{k}^{+})=%
\begin{cases}
1 & \text{if }l=\xi\\
0 & \text{if }l<\xi
\end{cases}
\\
&  \frac{\partial s_{j}}{\partial\omega_{jl}}(\tau_{k}^{+})=0\quad
\text{for all }l\leq\xi
\end{align}

\emph{Case 2}: $u_{j}(\tau_{k}^{-})=0,u_{j}(\tau_{k}^{+})=\pm1$
\begin{equation}
\frac{\partial s_{j}}{\partial\theta_{jl}}(\tau_{k}^{+})\hspace{-0.1cm}%
=\hspace{-0.1cm}%
\begin{cases}
\frac{\partial s_{j}}{\partial\theta_{jl}}(\tau_{k}^{-})-u_{j}(\tau_{k}%
^{+})sgn\big(\theta_{j\xi}-\theta_{j(\xi-1)}\big) & \hspace{-0.3cm} \text{if
}l=\xi\\
\frac{\partial s_{j}}{\partial\theta_{jl}}(\tau_{k}^{-})-u_{j}(\tau_{k}%
^{+})\Big[sgn(\theta_{jl}-\theta_{j(l-1)}) & \\
\quad\quad\quad\quad\quad\quad\quad-sgn(\theta_{j(l+1)}-\theta_{jl})\Big] &
\hspace{-0.3cm} \text{if }l<\xi\\
&
\end{cases}
\end{equation}%
\begin{equation}
\frac{\partial s_{j}}{\partial\omega_{jl}}(\tau_{k}^{+})=-u_{j}(\tau_{k}%
^{+})\quad\text{for all }l\leq\xi
\end{equation}

\emph{Case 3}: $u_{j}(\tau_{k}^{-})=\pm1,u_{j}(\tau_{k}^{+})=\mp1$
\begin{equation}
\frac{\partial s_{j}}{\partial\theta_{jl}}(\tau_{k}^{+})=%
\begin{cases}
2 & \text{if }l=\xi\\
-\frac{\partial s_{j}}{\partial\theta_{jl}}(\tau_{k}^{-}) & \text{if
}l<\xi
\end{cases}
\end{equation}

Details of these derivations can be found in \cite{cassandras2013optimal}. An
important difference arises in \emph{Case 2} above, where $\tau_{k}%
=|\theta_{j1}-a|+\omega_{j1}+...+|\theta_{j\xi}-\theta_{j(\xi-1)}%
|+\omega_{j\xi}$. We eliminate the constraints on the switching location that
$\theta_{j\xi}\leq\theta_{j(\xi-1)}$ if $\xi$ is even and $\theta_{j\xi}%
\geq\theta_{j(\xi-1)}$ if $\xi$ is odd.

\begin{figure}[h]
\begin{center}
\includegraphics[width=0.5\textwidth]{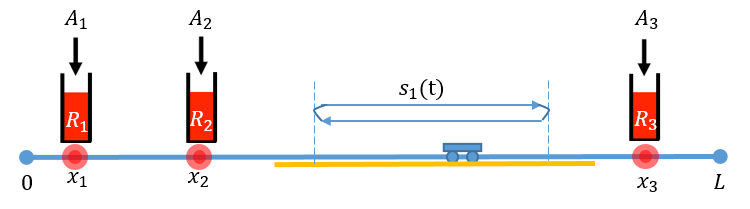}
\end{center}
\caption{\small An example of no event excitation leading to a failure of IPA
finding an optimal agent trajectory. The yellow bar is the segment of the
space covered by the agent.}%
\label{fig:1DimIPAFail}%
\end{figure}

\textbf{The event excitation problem}. Note that all derivative updates above
are limited to events occurring at times $\tau_{k}({\bm\theta},\bm\omega)$,
$k=1,2,\ldots$. Thus, this approach is scalable in the number of events
characterizing the hybrid system, not its state space. While this is a
distinct advantage, it also involves a potential drawback. In particular, it
assumes that the events involved in IPA updates are observable along a state
trajectory. However, if the current trajectory never reaches the vicinity of
any target so as to be able to sense it and affect the overall uncertainty
cost function, then any small perturbation to the trajectory will have no
effect on the cost. As a result, IPA will fail as illustrated in Fig.
\ref{fig:1DimIPAFail}: here, the single agent trajectory $s_{1}({\bm\theta
},\bm\omega,t)$ is limited to include no event. Thus, if a gradient-based
procedure is initialized with such $s_{1}({\bm\theta},\bm\omega,t)$, no event
involved in the evaluation of $\nabla R_{i}(t)$ is \textquotedblleft
excited\textquotedblright\ and the cost gradient remains zero.

In order to overcome this problem, we propose a modification of our cost
metric by introducing a function $V(\cdot)$ with the property of
\textquotedblleft spreading\textquotedblright\ the value of some $R_{i}(t)$
over all points $w\in\Omega\equiv\lbrack0,L]$. Recalling Proposition 1, we
limit ourselves to the subset $\mathcal{B}=[x_{1},x_{M}]\subset\Omega$. Then,
for all points $w\in\mathcal{B}$, we define $V(w,t)$ as a continuous density
function which results in a total value equivalent to the weighted sum of the
target values $\sum_{i=1}^{M}R_{i}(t)$. We impose the condition that $V(w,t)$
be monotonically decreasing in the Euclidean distance $\Vert w-x_{i}\Vert$.
More precisely, we define $d_{i}^{+}(w)=\max\big(\Vert w-x_{i}\Vert,r\big)$
where $r=\min_{j=1,\ldots,N}\{r_{j}\}$ which ensures that $d_{i}^{+}(w) \geq r$.
Thus, $d_{i}^{+}(w)=r>0$ is fixed for all points within the target's vicinity, $w\in\lbrack x_{i}-r,x_{i}+r]$. We define%
\begin{equation}
V(w,t)=\sum_{i=1}^{M}\frac{\alpha_{i}R_{i}(t)}{d_{i}^{+}(w)}\label{eq:J2R}%
\end{equation}
Note that $V(w,t)$ corresponds to the \textquotedblleft total weighted reward
density\textquotedblright\ at $w\in\mathcal{B}$. The weight $\alpha_{i}$ may
be included to capture the relative importance of targets, but we shall
henceforth set $\alpha_{i}=1$ for all $i=1,\ldots,M$ for simplicity. In order
to differentiate points $w\in\mathcal{B}$ in terms of their location relative
to the agents states $s_{j}(t)$, $j=1,\ldots,N$, we also define the travel cost function
\begin{equation}
Q(w,\mathbf{s}(t))=\sum_{j=1}^{N}\Vert s_{j}(t)-w\Vert\label{eq:Pfunction2}%
\end{equation}
Using these definitions we introduce a new objective function component, which
is added to the objective function in (\ref{eq:costfunction}):
\begin{equation}
J_{2}(t)=\int_{\mathcal{B}}Q(w,\mathbf{s}(t))V(w,t)dw\label{eq:J2}%
\end{equation}
The significance of $J_{2}(t)$ is that it accounts for the movement of agents
through $Q(w,{\mathbf{s}}(t))$ and captures the target state values through
$V(w,t)$. Introducing this term in the objective function in the following creates a
non-zero gradient even if the agent trajectories are not passing through any
targets. We now define the metrics in (\ref{eq:ParametricObj}) as $J_{1}(t)$ and incorporate the parametric $J_{2}(t)$ as an addition.
\begin{equation}
\hspace{-0.1cm}\min\limits_{\bm\theta\in\Theta,\bm\omega\geq0}\hspace
{-0.3cm}J(\bm\theta,\bm\omega,T)=\frac{1}{T}\int_{0}^{T}\big[J_{1}%
(\bm\theta,\bm\omega,t)+e^{-\beta t}J_{2}(\bm\theta,\bm\omega,t)\big]dt%
\label{eq:ParamOptim2}%
\end{equation}
where $J_{1}(\bm\theta,\bm\omega,t)=\sum_{i=1}^{M}R_{i}(t)$ is the original
uncertainty metric. This creates a continuous potential field for the agents
which ensures a non-zero cost gradient even when the trajectories do not
excite any events. This non-zero gradient will induce trajectory adjustments
that naturally bring them toward ones with observable events. The inclusion of
the factor $e^{-\beta t}$ with $\beta>0$ is included so that as the number of
IPA iterations increases, the effect of $J_{2}(\bm\theta,\bm\omega,t)$ is
diminished and the original objective is ultimately recovered. The IPA
derivative of $J_{2}(\bm\theta,\bm\omega,t)$ is
\begin{align}
&  \frac{\partial}{\partial\bm\theta}\int_{\tau_{k}}^{\tau_{k+1}}%
\int_{\mathcal{B}}Q(w,\bm\theta,\bm\omega,\mathbf{s}(t),t)V(w,\bm\theta
,\bm\omega,t)dw\nonumber\\%
&  =\int_{\tau_{k}}^{\tau_{k+1}}\int_{\mathcal{B}}\Big[\frac{\partial
Q(w,\bm\theta,\bm\omega,\mathbf{s}(t),t)}{\partial\bm\theta}V(w,\bm\theta,\bm\omega,t)\\%
& +Q(w,\bm\theta,\bm\omega,\mathbf{s}(t),t)\frac{\partial V(w,\bm\theta,\bm\omega
,t)}{\partial\bm\theta}\Big]dw\label{eq:IPAJ2der}%
\end{align}
where the derivatives of $Q(w,\bm\theta,\bm\omega,\mathbf{s}(t),t)$ and
$V(w,\bm\theta,\bm\omega,t)$ are obtained following the same procedure
described previously. Before making this modification, the lack of event
excitation in a state trajectory results in the total derivative
(\ref{eq:GradientParametricObj}) being zero. On the other hand, in
(\ref{eq:IPAJ2der}) we observe that if no events occur, the second part in the
integral, which involves $\frac{\partial V(\cdot)}{\partial\bm\theta}$ is
zero, since $\sum_{i=1}^{M}\frac{\partial R_{i}(t)}{\partial\bm\theta}=0$ at
all $t$. However, the first part in the integral does not depend on events,
but only the sensitivity of $Q(w,\bm\theta,\bm\omega,\mathbf{s}(t),t)$ in
(\ref{eq:Pfunction2}) with respect to the parameters $\bm\theta,\bm\omega$. As a
result, agent trajectories are adjusted so as to eventually excite desired
events and any gradient-based procedure we use in conjunction with IPA is no
longer limited by the absence of event excitation.

\textbf{IPA robustness to uncertainty modeling}. Observe that the evaluation
of $\nabla R_{i}\left(  t\right)  $, hence $\nabla J(\bm\theta,\bm\omega)$, is
\emph{independent} of $A_{i}$, $i=1,\ldots,M$, i.e., the parameters in our
uncertainty model. In fact, the dependence of $\nabla R_{i}\left(  t\right)  $
on $A_{i}$, $i=1,\ldots,M$, manifests itself through the event times $\tau
_{k}$, $k=1,\ldots,K$, that do affect this evaluation, but they, unlike
$A_{i}$ which may be unknown, are directly observable during the gradient
evaluation process. Thus, the IPA approach possesses an inherent
\emph{robustness} property: there is no need to explicitly model how
uncertainty affects $R_{i}(t)$ in (\ref{eq:multiDynR}). Consequently, we may
treat $A_{i}$ as unknown without affecting the solution approach (the values
of $\nabla R_{i}\left(  t\right)  $ are obviously affected). We may also allow
this uncertainty to be modeled through random processes $\{A_{i}(t)\}$,
$i=1,\ldots,M$. Under mild technical conditions on the statistical
characteristics of $\{A_{i}(t)\}$\cite{cassandras2010perturbation}, the
resulting $\nabla J(\bm\theta,\bm\omega)$ is an unbiased estimate of a
stochastic gradient.

\section{Graph-based scheduling method}
\label{sec:graphMethod}

While the IPA-driven gradient-based approach described in Sec. \ref{sec:IPA} offers several compelling advantages, it is not guaranteed to find a global optimum. In addition, it has been shown that in mission spaces of dimension greater than one, optimal trajectories cannot be described parametrically \cite{lin2013optimal}. This motivates the use of an alternative approach where the targets are viewed as discrete tasks, leading naturally to a graph-based description of the problem \cite{lan2013planning,smith2011persistent,lahijanian2010motion,smith2011optimal,mathew2013graph}. This higher level of abstraction allows one to guarantee an optimal solution, though at the cost of a significant increase in computational complexity. It is worth highlighting, however, that the complexity of such schemes is driven by the size of the graph and they are thus essentially invariant to the underlying dimensionality of the mission space.

\begin{figure}[htbp]
\centering \includegraphics[width=3.0in]{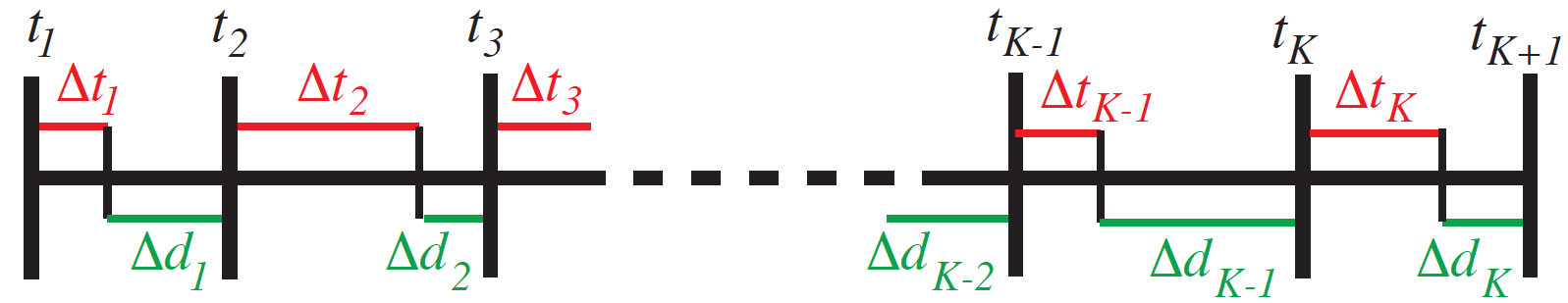}
\caption{\small Time sequence of a single agent on a given trajectory. The $t_i$ are the time points where the agent begins to move to the next target in the sequence. Each move takes $\Delta t_i$ units of time followed by a dwell period of $\Delta d_i$ units of time during which information is collected from the target.}
\label{fig:timeSequence}
\end{figure}
As illustrated in Fig. \ref{fig:timeSequence}, our approach to the discrete setting is to divide the overall planning time horizon $T$ for agent $j$ into a sum of $K_j$ consecutive time steps $\{t_j^1,t_j^2,...,t_j^{K_j}\}$, $j = 1,\dots,N$, with $t_j^1 = 0$. The dependence on $j$ indicates that each agent may have a different discretization.  We denote the end of the $K$-th step as $t_j^{K+1}=T$. Each step $k\in \{1,...,K_j\}$ begins with a travel stage where the agent moves to a particular target $i$. Under the assumption that during the transition between targets each agent moves at its maximum speed of $|u_j|=1$, the travel time is 
\begin{align}
	\Delta t_j^k = |s_j^k(t_j^k) - x_i|.\label{eq:TravelTime}
\end{align}
Upon arriving at a target, the agent dwells for a time $\Delta d_j^k$. Note that due to the range-based nature of the sensing, the uncertainty actually begins to decrease before the arrival of the agent at the target and continues to decrease after the agent has departed until the target is out of the sensing range.

The problem of optimizing the $u_j$ to minimize the average uncertainty over all the targets has been translated into a mixed integer programming (MIP) problem to select the sequence of targets and the dwell time at each target. Letting $a_{ji}^k$ be a binary variable denoting whether agent $j$ is assigned to target $i$ at time step $k$, this MIP is
\begin{align}
\min_{a_{ji}^k, \Delta d_j^k}& J = \frac{1}{T}\sum_{i=1}^M \int_0^T R_i(t)dt \label{eq:ObjDiscrete}\\
\textit{s.t.}& \quad a_{ji}^k \in \{1,0\}, \quad \sum_{i=1}^M a_{ji}^k = 1, \quad \forall j, k \\
&\quad \sum_{k=1}^K \Delta t_j^k + \Delta d_j^k \leq T,  \quad \forall j. \label{eq:TimeBound}
\end{align}

Note that we assume that each agent is assigned to a maximum of only one target at any one time. The IPA-driven approach has no such restriction. We break the solution of this problem into three parts: enumeration of all feasible trajectories, calculation of the cost of the feasible trajectories, and then selection of the optimal trajectory based on those costs. We focus on the case of a single agent for simplicity of description before generalizing to the multiple agent case. 

The first part, namely determining feasible trajectories, is straightforward. Given the fixed time horizon $T$, the target locations, the locations of the agent at the start of the time horizon, and the maximum speed of the agent,  a feasible trajectory is one where the sequence of targets can all be visited within the time horizon. Similarly, the third part simply involves comparing the trajectories and selecting the one with the minimal cost.

In the second part, the cost of each feasible trajectory must be determined. Suppose we have a given feasible trajectory with $K$ targets in its sequence. Note that because a trajectory may include multiple visits to the same target, $K$ may be larger than $m$ (and may be much larger for large time horizons and small $m$). Let $\{i_1, i_2, \dots, i_K \}$ denote the indices of the targets in the sequence. From \eqref{eq:ObjDiscrete}, the cost of this trajectory is given by the optimization problem
\begin{align*}
\min_{\Delta d_j^k}& J = \frac{1}{T}\sum_{i=1}^M \int_0^T R_i(t)dt \\
\textit{s.t.} &\quad \sum_{k=1}^K \Delta t^k + \Delta d^k \leq T.
\end{align*}
Our approach to solving this optimization problem is to setup a recursive calculation. As illustrated in Fig. \ref{fig:timeSequence}, since the travel times $\Delta t_i$ are completely determined by the sequence alone, optimizing over the dwell times is equivalent to optimizing the switching times $t_i$. Assume for the moment that the switching times through $t_{K-1}$ have been determined (and thus the first $K-2$ dwell times, $\Delta d^1, \dots, \Delta d^{K-2}$ are known). The two final dwell times are completely determined by selecting the time $t_K$ at which to switch the agent from target $i_{K-1}$ to target $i_K$. This then gives us a simple single variable optimization problem
\begin{align*}
\min_{\Delta T_K}& J = \frac{1}{\Delta T} \int_{t^{K-1}}^T (R_{i_{K-1}}(t)+R_{i_{K}}(t))\,dt
\end{align*}
where $\Delta T = T - t_{K-1}$. This allows the final switching time to be expressed as a function of the previous time $t_K = t_K(t_{K-1})$. Repeating this leads to an expression of the optimal switching times as a nested sequence of optimization functions which can be solved numerically.

This same optimization procedure can be generalized to the case of multiple agents. The primary challenge is that the set of feasible trajectories, and the calculation of the cost of those trajectories, quickly becomes intractable since all possible combinations of assignments of multiple agents must be considered. The computational complexity can be mitigated somewhat by taking advantage of known properties of optimal solutions (as described in Sec. \ref{sec:optimal}). 

Since the computationally complexity is exponential in the length of the time horizon, this approach is most feasible over short horizons. In prior work on linear systems, it was shown that an appropriately defined periodic schedule is sufficient to ensure the entire system remains controllable  \cite{yu2013effect,yu2014preservation}. In the current context, this translates to being able to keep the uncertainty of each of the targets arbitrarily close to zero. Motivated by this, we typically apply our discrete approach over a relatively short time horizon. If the resulting optimal trajectory is periodic, we extend it to longer horizons by simply repeating it.

\section{Simulation Examples}

To demonstrate the performance of the gradient-based algorithm using the IPA scheme described in Sec. \ref{sec:IPA}, we present two sets of numerical examples. The first set uses deterministic target locations and dynamics. The results are compared against the optimal found by the discrete scheduling algorithm of Sec. \ref{sec:graphMethod}. The second set demonstrates the robustness of the IPA scheme with respect to a stochastic uncertainty model.

\begin{figure}[htpb]
\centering
\includegraphics[width=\columnwidth]{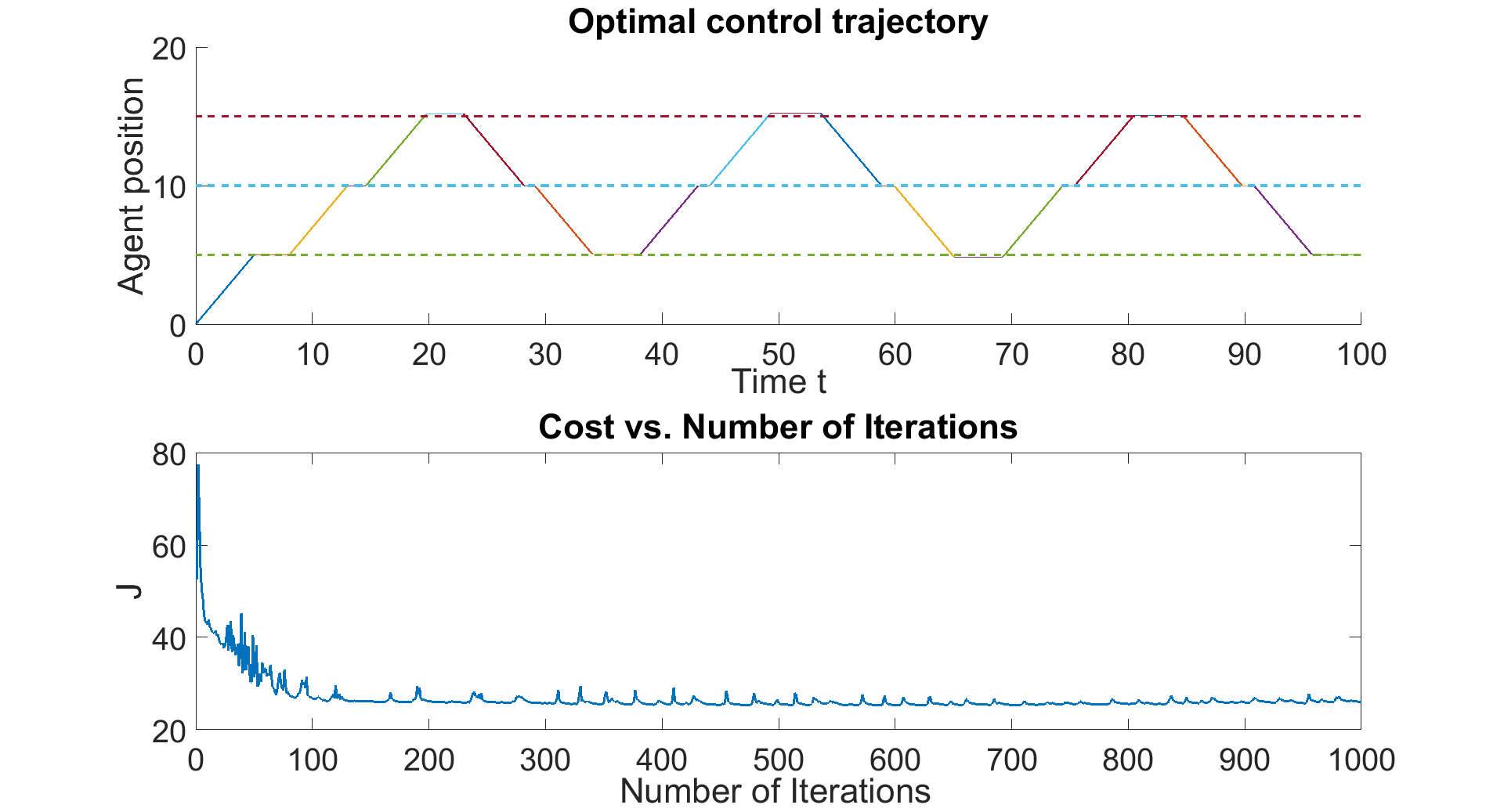} \\
\includegraphics[width=\columnwidth]{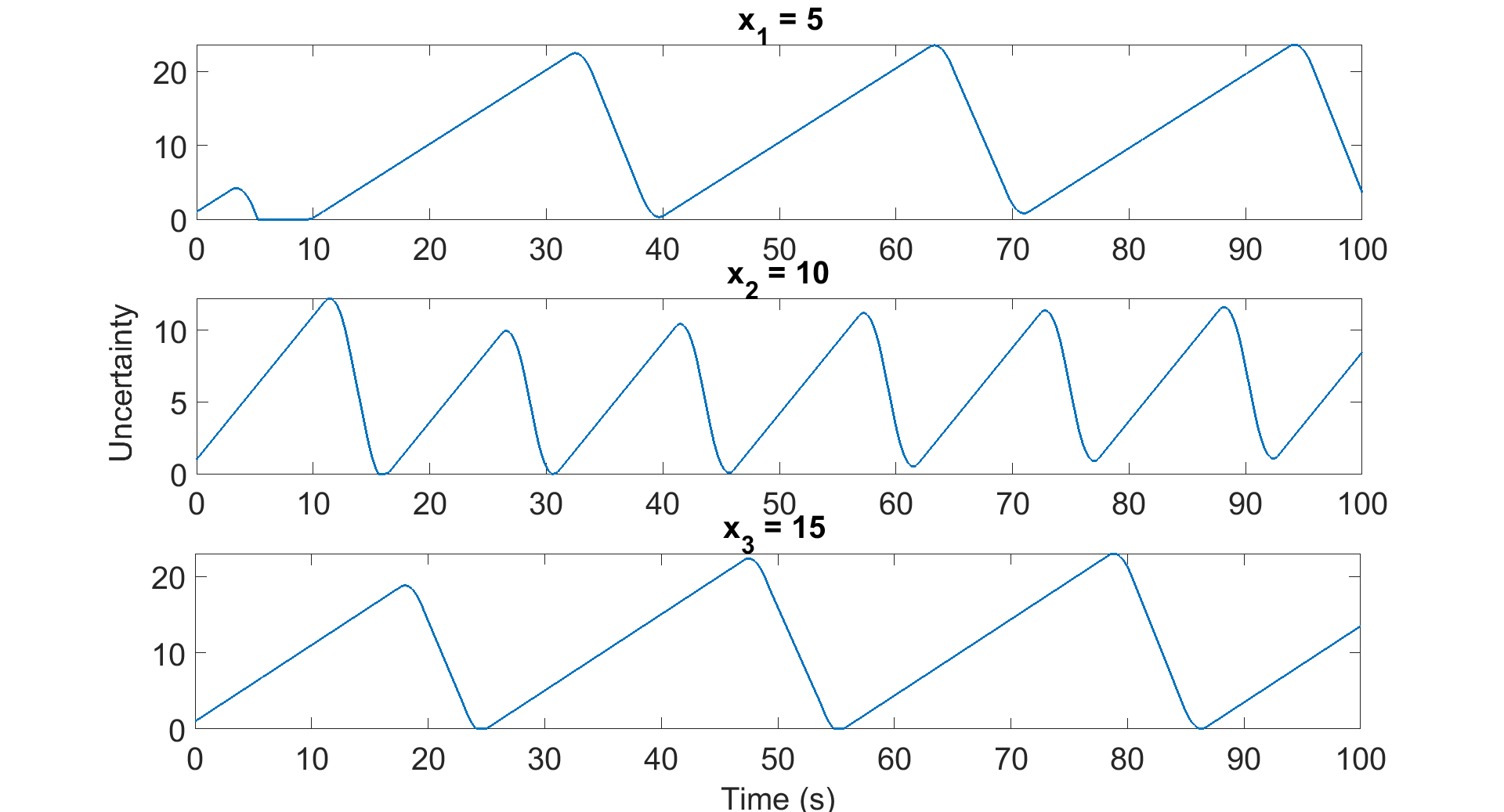}
\caption{\small  A single agent monitoring three targets using the IPA-driven gradient descent algorithm. (top image) Agent trajectory. (second image) Calculated cost as a function of iteration in the gradient descent. The final cost is 26.11. (bottom images) Target uncertainties along the trajectory.}
\label{fig:1A3T_IPA}
\end{figure}

The first simulation consists of a single agent performing a persistent monitoring task on three targets over a time horizon of 100 seconds. The targets are located at positions $x_{1}=5,$ $x_{2}=10,$ $x_{3}=15$ and their uncertainty dynamics in \eqref{eq:multiDynR} are defined by the parameters $A_i=1$, $B_i=5$, and $R_i(0) = 1$ for all $i$. The agent has a sensing range of 2 and is initialized with $s(0) = 0$, $u(0)=1$. The results from the IPA gradient descent approach are shown in Fig. \ref{fig:1A3T_IPA}. The top image shows the optimal trajectory of the agent determined after 1000 iterations of the IPA gradient descent while the bottom shows the evolution of the overall cost as a function of iteration number. The agent is moving through a periodic cycle of $x_1\rightarrow x_2 \rightarrow x_3 \rightarrow x_2 \rightarrow x_1\cdots$, dwelling for a short time at each target before moving to the next. Notice that the agent dwells for a shorter time at the center target since it visits that location twice per cycle. The second image in the figure shows that gradient descent converges within the first 100 iterations. This simulation aims to test the event driven IPA scheme with the discrete scheduling algorithm which yields optimal but suffers from computational intensity. Thus, we start with a short time horizon $T=100$s. Event-driven IPA optimizes the trajectory fast but the convergence is somewhat unstable due to the lack of events within a short time horizon. The final cost is 26.11. The bottom images in Fig. \ref{fig:1A3T_IPA} show the evolution of the target uncertainties. 

\begin{figure}[htpb]
\centering
\includegraphics[width=\columnwidth]{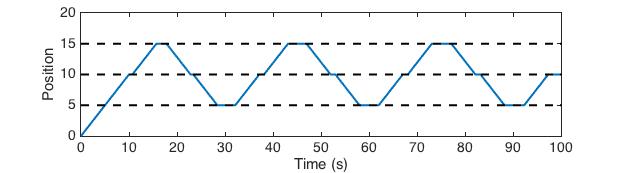} \\
\includegraphics[width=\columnwidth]{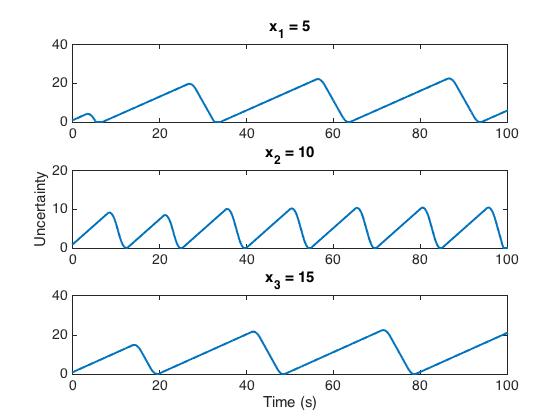} \\
\caption{\small A single agent monitoring three targets using the optimal discrete assignment and dwelling time. The final cost is 25.07. (top image) The agent trajectory is almost the same as in Fig. \ref{fig:1A3T_IPA}. (bottom images) Target uncertainties along the trajectory.}
\label{fig:1A3T_Graph}
\end{figure}

The corresponding result based on the discrete setting of Sec.
\ref{sec:graphMethod} is essentially the same with the agent moving through
the three targets in a periodic fashion as shown in Fig. \ref{fig:1A3T_Graph}. The only deviation from the IPA scheme occurs at the end of the horizon where the discrete approach returns to the center target. The final cost was 25.07, matching that of the IPA approach and thus verifying the approximate optimality of the solution found in Fig. \ref{fig:1A3T_IPA}. 

\begin{figure}[htpb]
\centering
\includegraphics[width=\columnwidth]{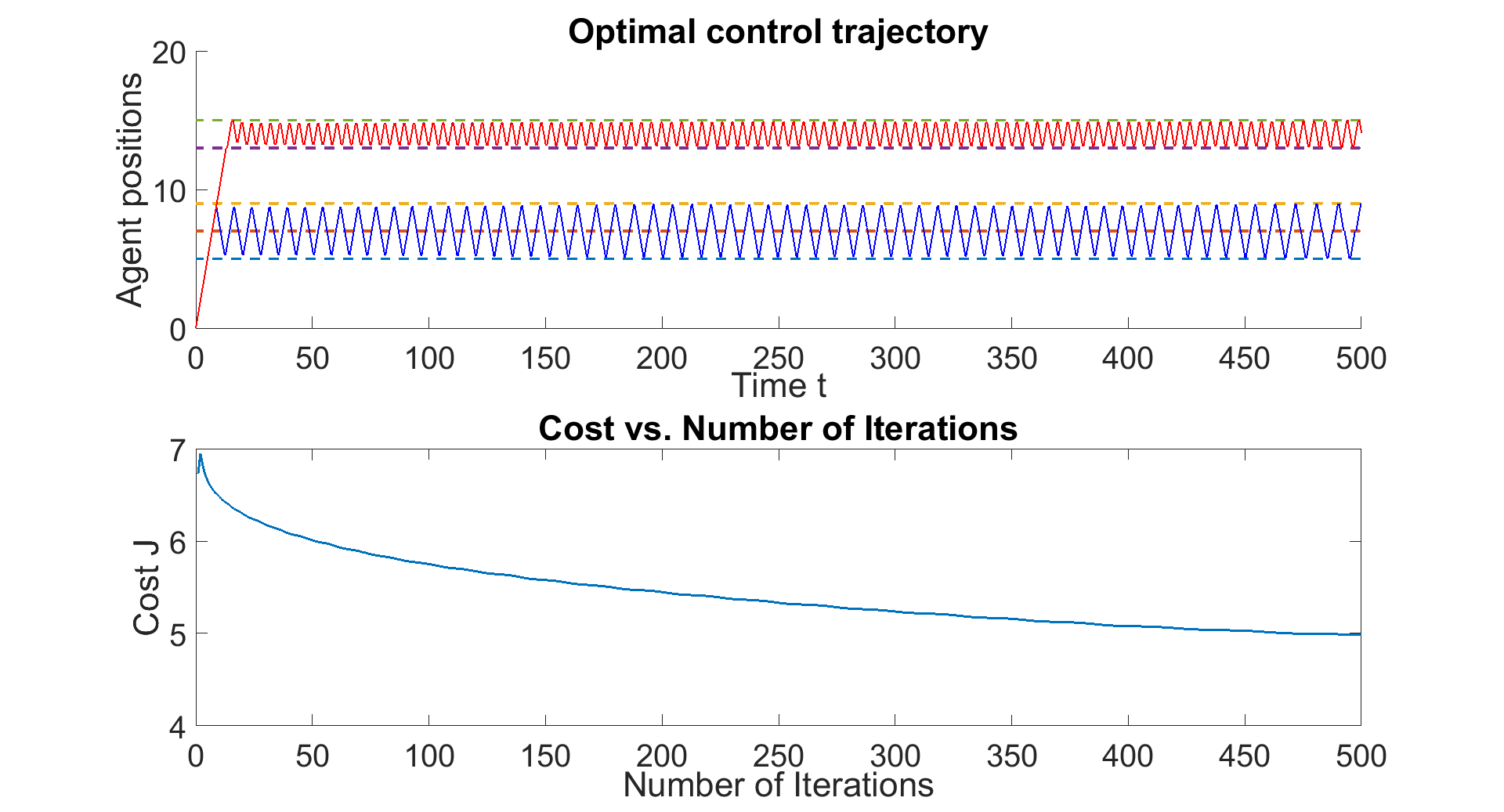} \\
\includegraphics[width=\columnwidth]{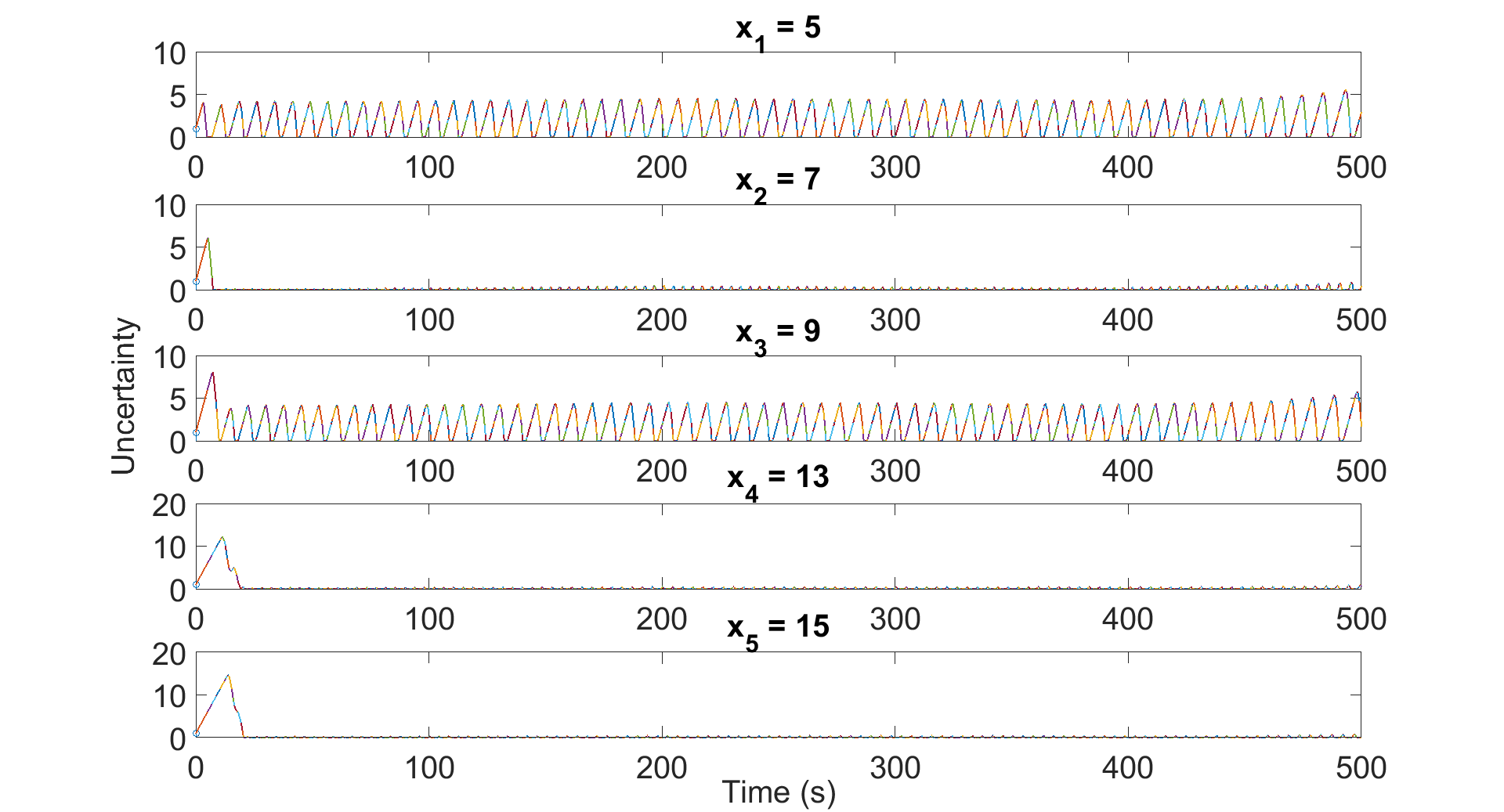}	
\caption{\small Two agents monitoring five targets using the IPA gradient descent algorithm. (top image) Agent trajectories. (second image) Calculated cost as a function of iteration. The final cost is 4.99. (bottom images) Target uncertainty values along the above trajectories.}
\label{fig:2A5T_IPA}
\end{figure}

\begin{figure}[htpb]
\centering
\includegraphics[width=\columnwidth]{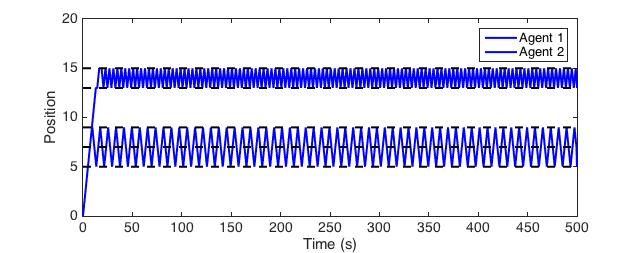} \\
\includegraphics[width=\columnwidth]{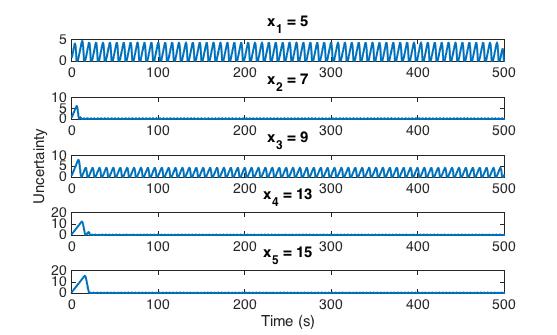}	
\caption{\small Two agents monitoring five targets using the discrete assignment and dwelling time. The final cost was 4.92. (top image) Agent trajectories. (bottom images) Target uncertainty values along the above trajectories.}
\label{fig:2A5T_graph}
\end{figure}

The next simulation involves two agents and five targets over a time horizon of 500 seconds. The targets are located at $x_{1}=5,$ $x_{2}=7,$ $x_{3}=9,$ $x_{4}=13,$ $x_{5}=15$. The uncertain dynamics were the same as in the single agent, three target case. As before, the agents have a sensing range of 2 and are initialized at  $s_{1}(0)=s_{2}(0)=0,$ with $u_{1}(0)=u_{2}(0)=1$. The results from the event-driven IPA gradient descent approach are shown in Fig. \ref{fig:2A5T_IPA}. The solution is again periodic with the agents dividing the targets into two groups. Notice that the single agent on targets $x_4$ and $x_5$ is able to keep the uncertainties very close to zero since the targets are quite close relative to the sensing range of the agent. The other agent is able to hold its middle target ($x_2$) close to zero since it is visited more often. The uncertainties of targets $x_1$ and $x_3$ rise and decrease to zero constantly. The corresponding result based on the discrete setting is shown in Fig. \ref{fig:2A5T_graph}. Rather than solve over the full horizon, the problem was solved over a 60 second horizon and then the periodic trajectory repeated to fill the 500 second horizon. The results are again very close to the event-driven IPA method.

Note that the optimal trajectories in both one and two agent examples are bounded between $[5,15]$ (positions of the first and last target), which
is consistent with Proposition 1.

\begin{figure}[ptb]
\centering
\begin{subfigure}[b]{0.475\textwidth}
\includegraphics[width=\textwidth]{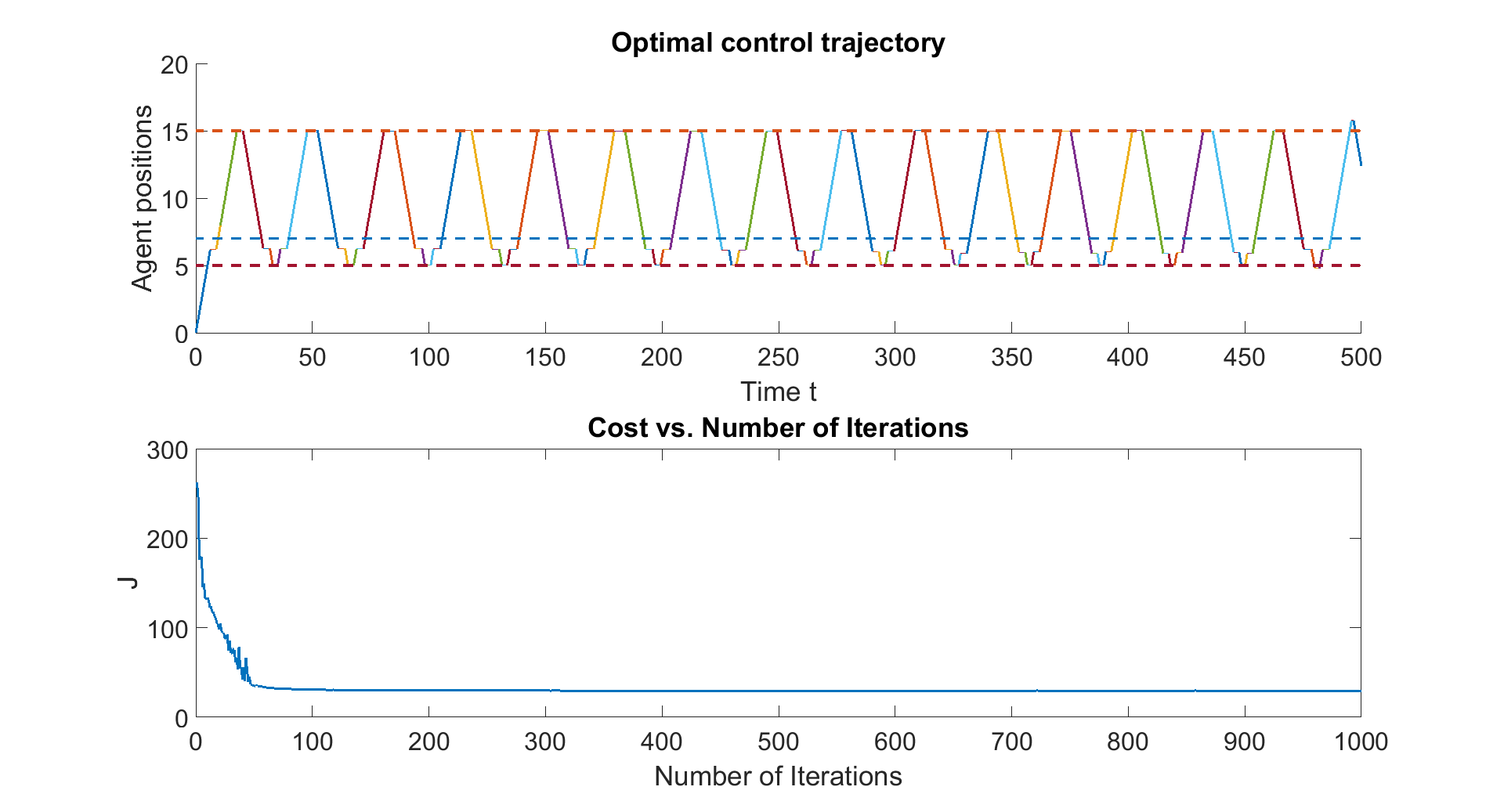}
\caption[]{{\small Example of deterministic target model. Target positions $5,7,15$, dynamics parameter $A_i = 1, B=5, r=2$. $J^*(\bm \theta,\bm \omega) = 29.40$.}}
\label{fig:Traj_1A3T}
\end{subfigure}
\centering
\begin{subfigure}[b]{0.475\textwidth}
\includegraphics[width=\textwidth]{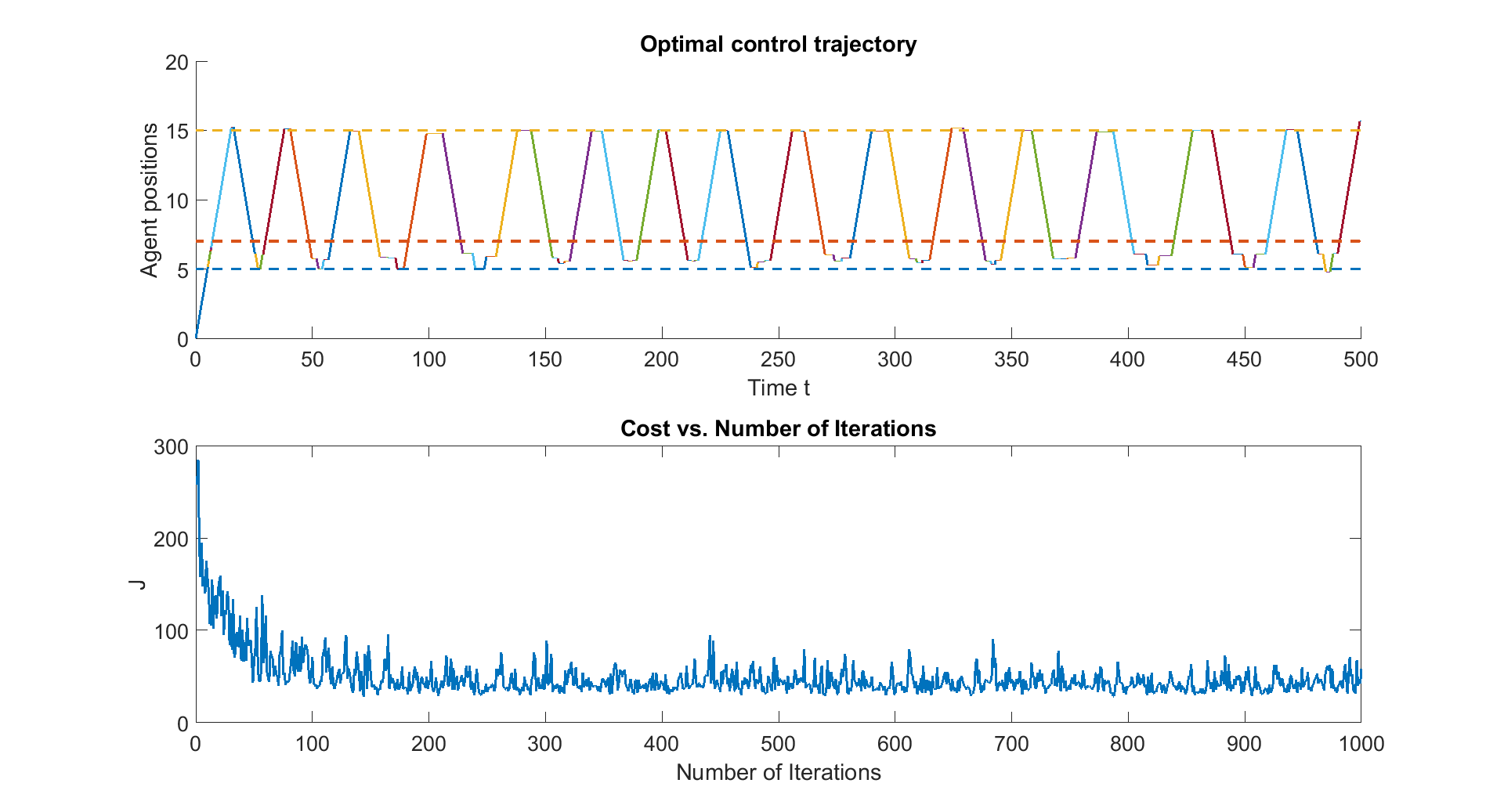}
\caption[]{\small Example with stochastic uncertainty processes. $A_i \sim U(0,2)$. $J^{\ast}(\bm\theta,\bm\omega)=42.46$.}
\label{fig:Inflow_stochastic_1A_traj}
\end{subfigure}
\centering
\begin{subfigure}[b]{0.475\textwidth}
\includegraphics[width=\textwidth]{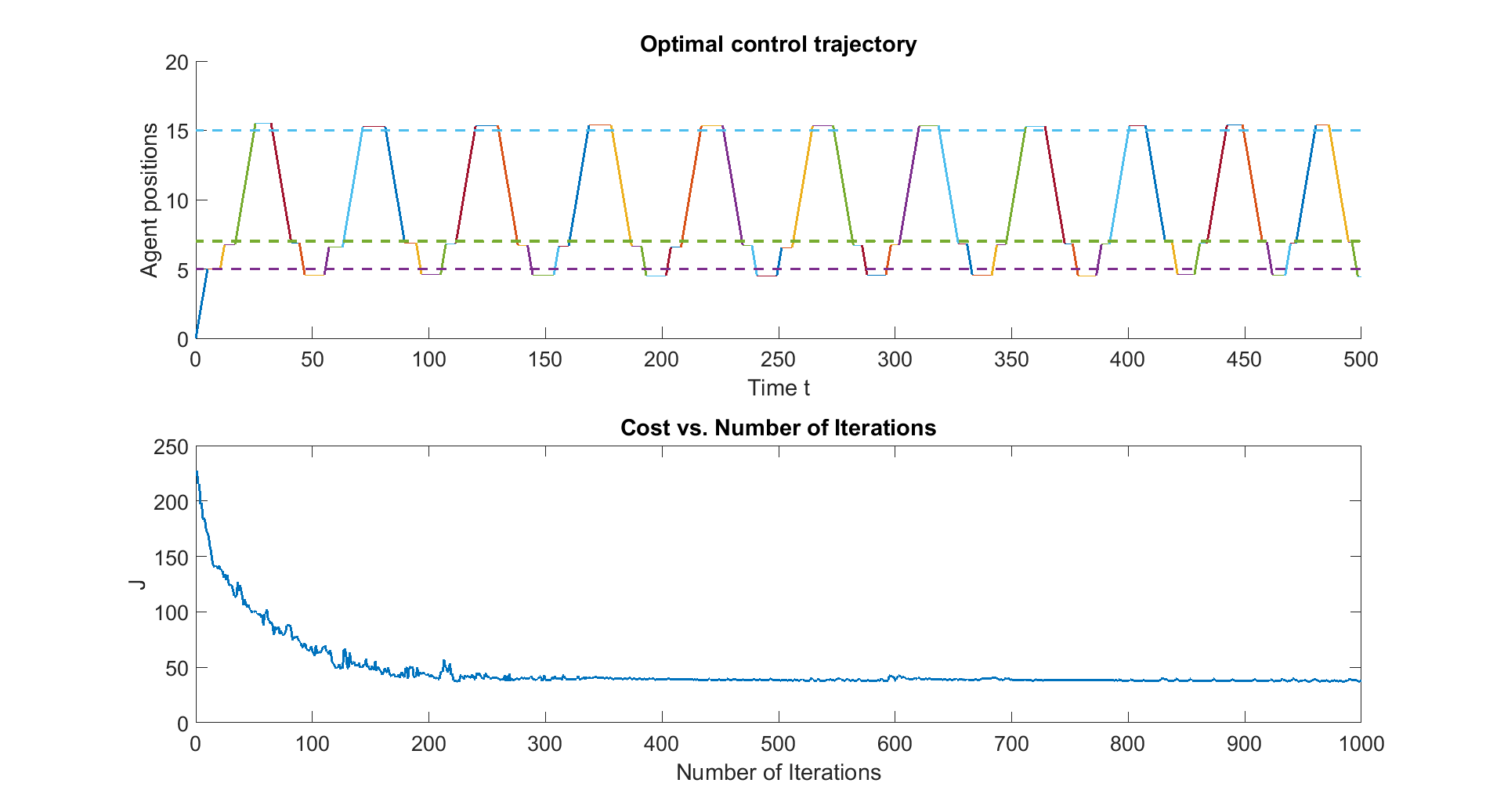}
\caption[]{\small Example with stochastic target locations $\sim U(x_i-0.25,x_i+0.25)$. $J^{\ast}(\bm\theta,\bm\omega)=34.89$.}
\label{fig:TarPos_stochastic_1A_traj}
\end{subfigure}
\caption[]{{\protect\small Examples demonstrating IPA robustness with respect to stochastic uncertainty. (a)(b)(c) Top plot: optimal trajectory $s^{\ast}(t)$. Bottom plot:
cost convergence. }}%
\label{fig:IPA_robustness}%
\end{figure}

As mentioned earlier, the IPA robustness property allows us to handle
stochastic uncertainty models at targets. We show a one-agent example in Fig.
\ref{fig:Inflow_stochastic_1A_traj} where the uncertainty inflow rate
$A_{i}(t)$ is uniformly distributed between $[0,2]$ for all targets. In Fig.
\ref{fig:TarPos_stochastic_1A_traj}, we introduce randomness by allowing target
positions to vary uniformly over $[x_{i}-0.25,x_{i}+0.25]$. In both cases, the
optimal cost in the stochastic models in Figs.
\ref{fig:Inflow_stochastic_1A_traj} and \ref{fig:TarPos_stochastic_1A_traj}
are close to the optimal cost of the deterministic case Fig.
\ref{fig:Traj_1A3T} where the parameter $A_{i}$ and target positions $x_{i}$
are the means of the associated random processes in the stochastic models. As
expected, the convergence depends on the variance of these random processes.

\begin{figure}[ptb]
\centering
\begin{subfigure}[b]{0.475\textwidth}
\includegraphics[width=\textwidth]{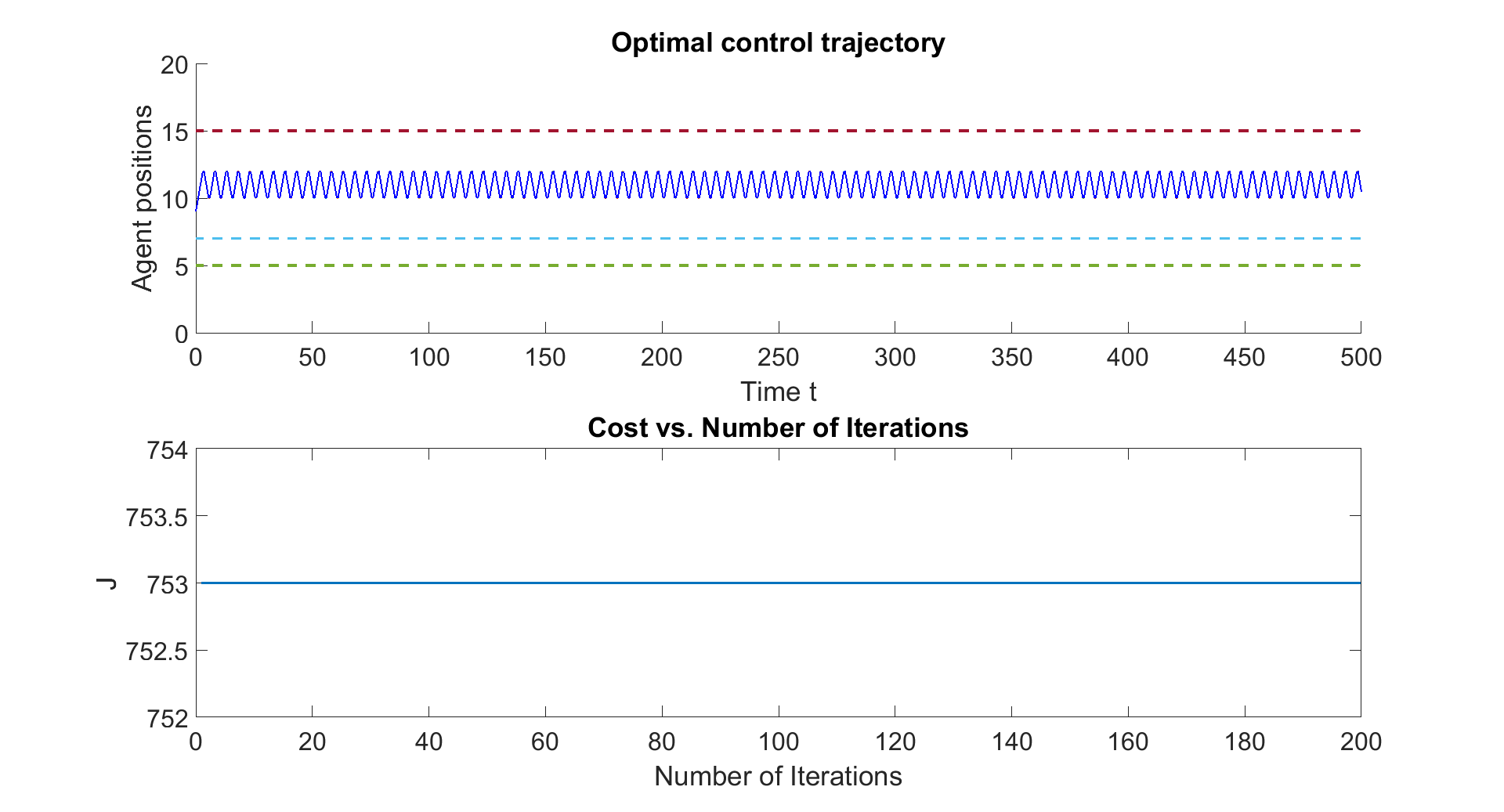}
\caption[]{{\small A trajectory where IPA fails due to lack of event excitation. Top plot: agent trajectory. Bottom plot: cost convergence.}}
\label{fig:IPAFail_1A_traj_png}
\end{subfigure}
\centering
\begin{subfigure}[b]{0.475\textwidth}
\includegraphics[width=\textwidth]{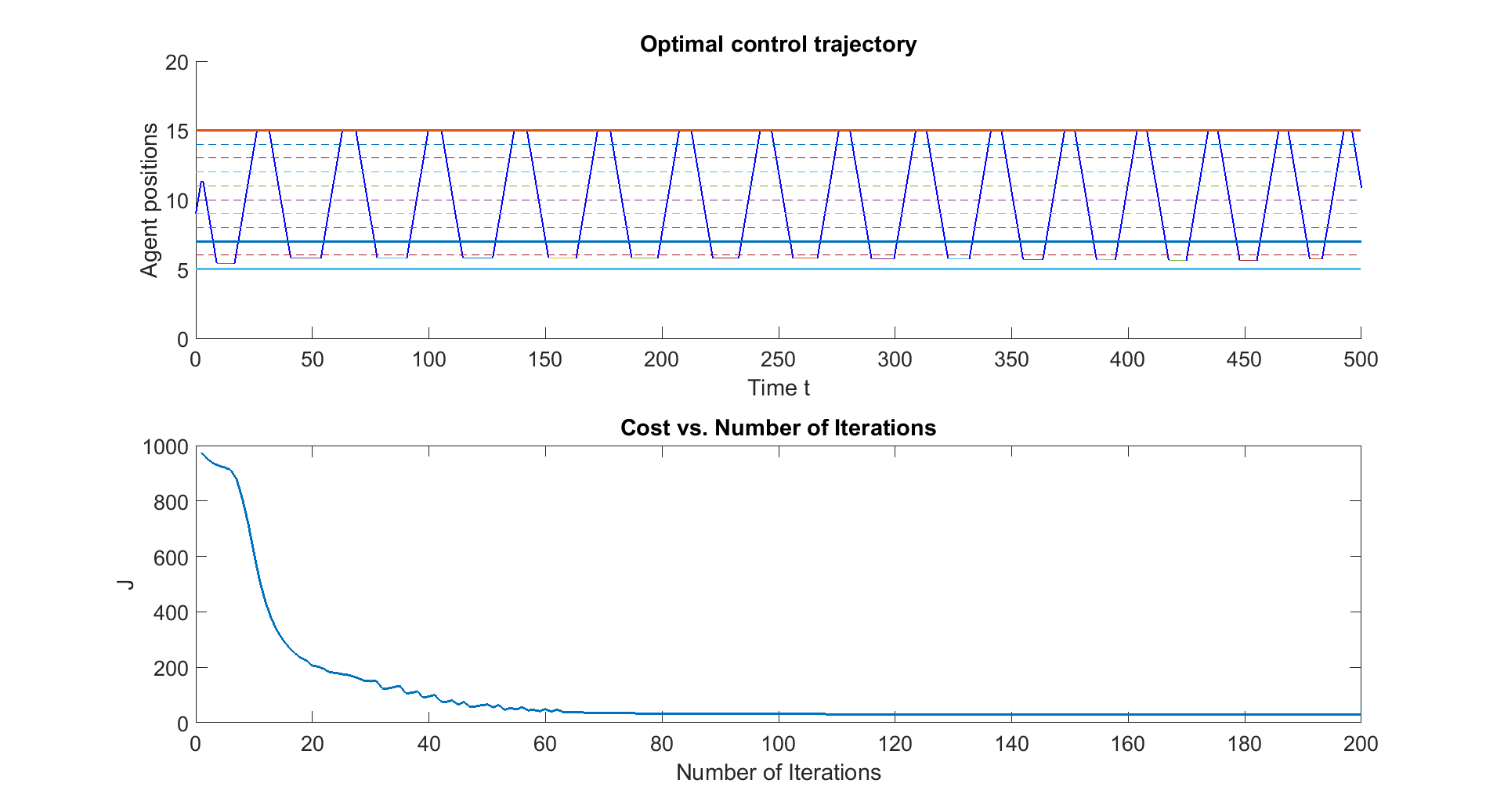}
\caption[]{{\small IPA optimization after event excitation. Top plot: optimal agent trajectory. Bottom plot: cost convergence. $J^{\ast}(\bm\theta,\bm\omega)=30.24$. }}
\label{fig:IPAFail_NewMetric_traj}
\end{subfigure}
\centering
\begin{subfigure}[b]{0.475\textwidth}
\includegraphics[width=\textwidth]{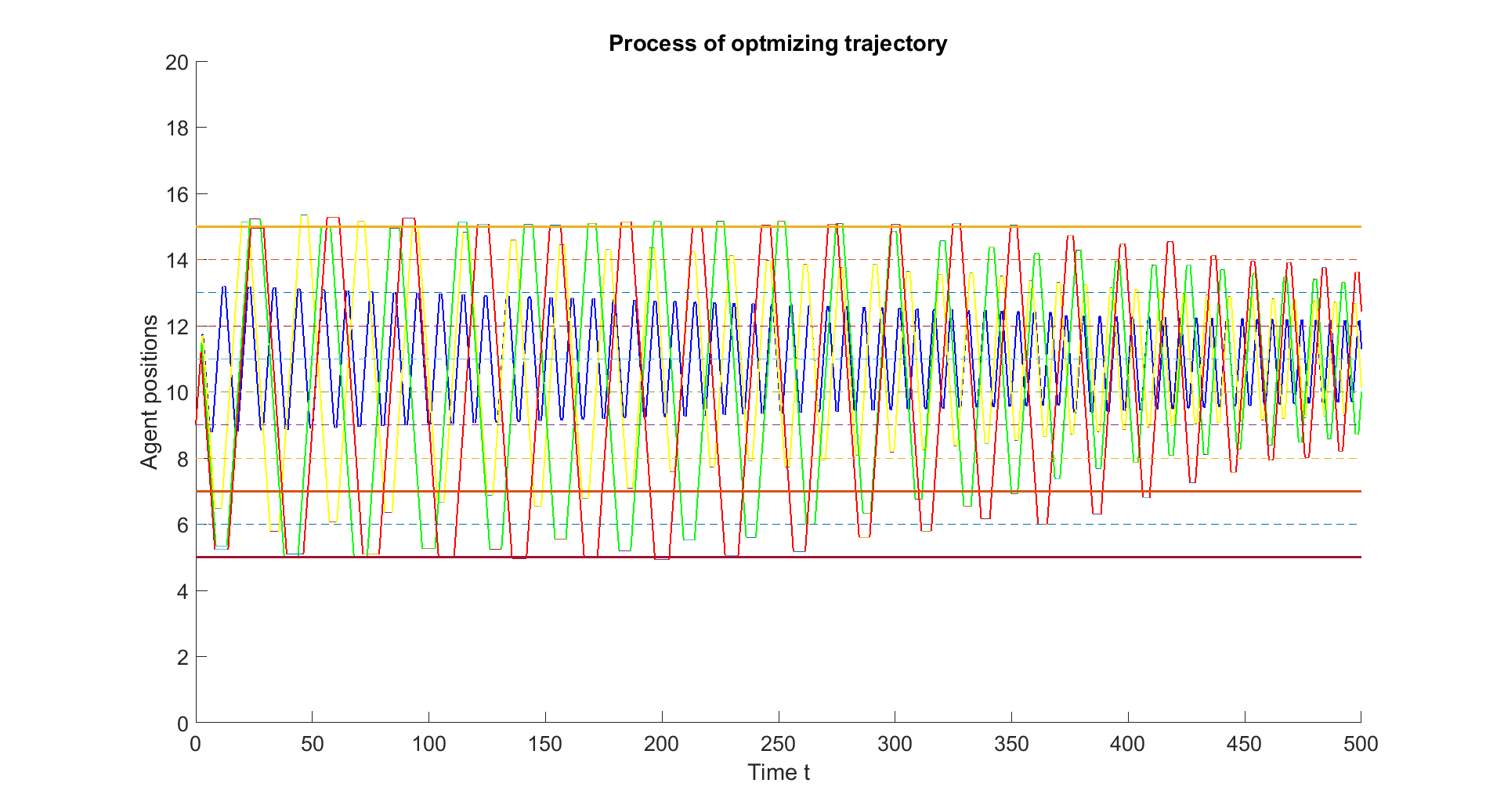}
\caption[]{{\small Trajectory adjustments with event excitation after 5 (blue), 10 (green), and 15 (red) iterations.}}
\label{fig:NewMetric_optimizing_process}
\end{subfigure}
\caption[]{{\protect\small The event excitation issue. After adding $J_{2}(\bm\theta,\bm\omega,t)$, the trajectory adjusts to include targets, the cost converges to 30.24 which is close to the optimal cost in Fig.\ref{fig:Traj_1A3T} where the target dynamics are the same. }}
\label{fig:IPAFail_NewMetric}%
\end{figure}

The event excitation issue is addressed in Fig. \ref{fig:IPAFail_1A_traj_png},
where the agent trajectory is initialized so that it is not close to any of
the targets. Using the original problem formulation (without the inclusion of
$J_{2}(\bm\theta,\bm\omega,t)$ in (\ref{eq:ParamOptim2}), the initial trajectory and cost remain unchanged. After adding $J_{2}(\bm\theta,\bm\omega,t)$, the
blue, green, and red curves in Fig. \ref{fig:NewMetric_optimizing_process}
show the trajectory adjustment after 5, 10, and 15 iterations respectively.
After 100 iterations, the cost converges to 30.24 as shown in Fig.
\ref{fig:IPAFail_NewMetric_traj} which is close to the optimal cost in Fig.
\ref{fig:Traj_1A3T} where the target dynamics are the same.

\section{Conclusion}

We have formulated a persistent monitoring problem with the objective of controlling the movement of multiple cooperating agents so as to minimize an uncertainty metric associated with a finite number of targets. We have established properties of the optimal control solution which reduce the problem to a parametric optimization one. A complete on-line solution is given
by Infinitesimal Perturbation Analysis (IPA) to evaluate the gradient of the
objective function with respect to all parameters. We also address the case when
IPA gradient estimation fails because of the lack of event excitation.  We solve this problem by proposing a new metric for the objective function which creates a potential field guaranteeing that gradient values are non-zero. This approach is
compared to an alternative graph-based task scheduling algorithm for
determining an optimal sequence of target visits. Ongoing research is investigating how to extend these methodologies to higher dimensional mission spaces.


\bibliographystyle{IEEEtran}
\bibliography{library}

\begin{thebibliography}{10}
\providecommand{\url}[1]{#1}
\csname url@samestyle\endcsname
\providecommand{\newblock}{\relax}
\providecommand{\bibinfo}[2]{#2}
\providecommand{\BIBentrySTDinterwordspacing}{\spaceskip=0pt\relax}
\providecommand{\BIBentryALTinterwordstretchfactor}{4}
\providecommand{\BIBentryALTinterwordspacing}{\spaceskip=\fontdimen2\font plus
\BIBentryALTinterwordstretchfactor\fontdimen3\font minus
  \fontdimen4\font\relax}
\providecommand{\BIBforeignlanguage}[2]{{%
\expandafter\ifx\csname l@#1\endcsname\relax
\typeout{** WARNING: IEEEtran.bst: No hyphenation pattern has been}%
\typeout{** loaded for the language `#1'. Using the pattern for}%
\typeout{** the default language instead.}%
\else
\language=\csname l@#1\endcsname
\fi
#2}}
\providecommand{\BIBdecl}{\relax}
\BIBdecl

\bibitem{zhong2011distributed}
M.~Zhong and C.~G. Cassandras, ``Distributed coverage control and data
  collection with mobile sensor networks,'' \emph{Automatic Control, IEEE
  Transactions on}, vol.~56, no.~10, pp. 2445--2455, 2011.

\bibitem{sun2015optimal}
X.~Sun and C.~G. Cassandras, ``Optimal dynamic formation control of multi-agent
  systems in environments with obstacles,'' \emph{arXiv preprint
  arXiv:1508.04727}, 2015.

\bibitem{cassandras2013optimal}
C.~Cassandras, X.~Lin, and X.~Ding, ``An optimal control approach to the
  multi-agent persistent monitoring problem,'' \emph{IEEE Transactions on
  Automatic Control}, vol.~58, no.~4, pp. 947--961, 2013.

\bibitem{lin2013optimal}
X.~Lin and C.~Cassandras, ``An optimal control approach to the multi-agent
  persistent monitoring problem in two-dimensional spaces,'' in \emph{Proc. of
  the IEEE Conference on Decision and Control}.\hskip 1em plus 0.5em minus
  0.4em\relax IEEE, 2013, pp. 6886--6891.

\bibitem{michael2011persistent}
N.~Michael, E.~Stump, and K.~Mohta, ``Persistent surveillance with a team of
  mavs,'' in \emph{2011 IEEE/RSJ International Conference on Intelligent Robots
  and Systems}, 2011.

\bibitem{smith2011persistent}
S.~L. Smith, M.~Schwager, and D.~Rus, ``Persistent monitoring of changing
  environments using a robot with limited range sensing,'' in \emph{Proc. of
  the IEEE International Conference on Robotics and Automation (ICRA)}.\hskip
  1em plus 0.5em minus 0.4em\relax IEEE, 2011, pp. 5448--5455.

\bibitem{Shen:2011kj}
Z.~Shen and S.~B. Andersson, ``{Tracking Nanometer-Scale Fluorescent Particles
  in Two Dimensions With a Confocal Microscope},'' \emph{IEEE Transactions on
  Control Systems Technology}, vol.~19, no.~5, pp. 1269--1278, Sep. 2011.

\bibitem{cromer2011tracking}
S.~M. Cromer~Berman, P.~Walczak, and J.~W. Bulte, ``Tracking stem cells using
  magnetic nanoparticles,'' \emph{Wiley Interdisciplinary Reviews: Nanomedicine
  and Nanobiotechnology}, vol.~3, no.~4, pp. 343--355, 2011.

\bibitem{horling2004survey}
B.~Horling and V.~Lesser, ``A survey of multi-agent organizational paradigms,''
  \emph{The Knowledge Engineering Review}, vol.~19, no.~04, pp. 281--316, 2004.

\bibitem{yu2015persistent}
J.~Yu, S.~Karaman, and D.~Rus, ``Persistent monitoring of events with
  stochastic arrivals at multiple stations,'' \emph{Robotics, IEEE Transactions
  on}, vol.~31, no.~3, pp. 521--535, 2015.

\bibitem{stump2011multi}
E.~Stump and N.~Michael, ``Multi-robot persistent surveillance planning as a
  vehicle routing problem,'' in \emph{Automation Science and Engineering
  (CASE), 2011 IEEE Conference on}.\hskip 1em plus 0.5em minus 0.4em\relax
  IEEE, 2011, pp. 569--575.

\bibitem{cassandras2010perturbation}
C.~G. Cassandras, Y.~Wardi, C.~G. Panayiotou, and C.~Yao, ``Perturbation
  analysis and optimization of stochastic hybrid systems,'' \emph{European
  Journal of Control}, vol.~16, no.~6, pp. 642--661, 2010.

\bibitem{wardi2010unified}
Y.~Wardi, R.~Adams, and B.~Melamed, ``A unified approach to infinitesimal
  perturbation analysis in stochastic flow models: the single-stage case,''
  \emph{Automatic Control, IEEE Transactions on}, vol.~55, no.~1, pp. 89--103,
  2010.

\bibitem{Schwager2009decentralized}
M.~Schwager, D.~Rus, and J.-J. Slotine, ``Decentralized, adaptive coverage
  control for networked robots,'' \emph{The International Journal of Robotics
  Research}, vol.~28, no.~3, pp. 357--375, 2009.

\bibitem{cao2011maintaining}
M.~Cao, A.~S. Morse, C.~Yu, B.~Anderson, S.~Dasgupta \emph{et~al.},
  ``Maintaining a directed, triangular formation of mobile autonomous agents,''
  \emph{Communications in Information and Systems}, vol.~11, no.~1, p.~1, 2011.

\bibitem{oh2014formation}
K.-K. Oh and H.-S. Ahn, ``Formation control and network localization via
  orientation alignment,'' \emph{Automatic Control, IEEE Transactions on},
  vol.~59, no.~2, pp. 540--545, 2014.

\bibitem{YasamanKhazaeni}
Y.~Khazaeni and C.~G. Cassandras, ``Event excitation for event-driven control
  and optimization of multi-agent systems,'' in \emph{IEEE International
  Workshop on Discrete Event Systems(WODES)}.\hskip 1em plus 0.5em minus
  0.4em\relax IEEE, 2016.

\bibitem{bryson1975applied}
A.~E. Bryson, \emph{Applied optimal control: optimization, estimation and
  control}.\hskip 1em plus 0.5em minus 0.4em\relax CRC Press, 1975.

\bibitem{lan2013planning}
X.~Lan and M.~Schwager, ``Planning periodic persistent monitoring trajectories
  for sensing robots in gaussian random fields,'' in \emph{IEEE International
  Conference on Robotics and Automation (ICRA)}.\hskip 1em plus 0.5em minus
  0.4em\relax IEEE, 2013, pp. 2415--2420.

\bibitem{lahijanian2010motion}
M.~Lahijanian, J.~Wasniewski, S.~B. Andersson, and C.~Belta, ``Motion planning
  and control from temporal logic specifications with probabilistic
  satisfaction guarantees,'' in \emph{IEEE International Conference on Robotics
  and Automation (ICRA)}.\hskip 1em plus 0.5em minus 0.4em\relax IEEE, 2010,
  pp. 3227--3232.

\bibitem{smith2011optimal}
S.~L. Smith, J.~Tumova, C.~Belta, and D.~Rus, ``Optimal path planning for
  surveillance with temporal logic constraints,'' \emph{The International
  Journal of Robotics Research}, p. 0278364911417911, 2011.

\bibitem{mathew2013graph}
N.~Mathew, S.~L. Smith, and S.~L. Waslander, ``A graph-based approach to
  multi-robot rendezvous for recharging in persistent tasks,'' in
  \emph{Robotics and Automation (ICRA), 2013 IEEE International Conference
  on}.\hskip 1em plus 0.5em minus 0.4em\relax IEEE, 2013, pp. 3497--3502.

\bibitem{yu2013effect}
X.~Yu and S.~B. Andersson, ``Effect of switching delay on a networked control
  system,'' in \emph{Proc.of the IEEE Conference on Decision and Control
  (CDC)}.\hskip 1em plus 0.5em minus 0.4em\relax IEEE, 2013, pp. 5945--5950.

\bibitem{yu2014preservation}
X.~Yu and S.~B. Andersson, ``Preservation of system properties for networked linear,
  time-invariant control systems in the presence of switching delays,'' in
  \emph{Proc. of the IEEE Conference on Decision and Control}.\hskip 1em plus
  0.5em minus 0.4em\relax IEEE, 2014, pp. 5260--5265.

\end{thebibliography}

\end{document}